\documentclass[12pt]{amsart}
\usepackage{latexsym, amssymb, amsmath}
 \usepackage{eucal}

    \newtheorem{rema}{Remark}[section]
    \newtheorem{propo}[rema]{Proposition}

   \newtheorem{theo}[rema]{Theorem}
 \newtheorem{conj}[rema]{Conjecture}

   \newtheorem{defi}[rema]{Definition}
    \newtheorem{lemma}[rema]{Lemma}
    \newtheorem{corol}[rema]{Corollary}
     
  \newtheorem{rmk}[rema]{Remark}

	\newcommand{\nno}{\nonumber}

	\newcommand{\p}{\partial}
\newcommand{\fr}{\frac}

 \newcommand{\pf}{{\it Proof:}\hspace{2ex}}
 
 \newcommand{\epfv}{\hspace{1em}$\Box$\vspace{1em}}

\newcommand{\bC}{{\mathbb C}}

\newcommand{\bQ}{{\mathbb Q}}
\newcommand{\bN}{{\mathbb N}}

\newcommand{\BQ}{\begin{eqnarray}}
\newcommand{\EQ}{\end{eqnarray}}
\newcommand{\BQn}{\begin{eqnarray*}}
\newcommand{\EQn}{\end{eqnarray*}}
\newcommand{\BL}{\begin{align}}
\newcommand{\EL}{\end{align}}
\newcommand{\BLn}{\begin{align*}}
\newcommand{\ELn}{\end{align*}}
\newcommand{\BA}{\begin{align}}
\newcommand{\EA}{\end{align}}
\newcommand{\BAn}{\begin{align*}}
\newcommand{\EAn}{\end{align*}}
\newcommand{\wtilde}{\widetilde}

\newcommand{\Hes}{ \text{Hes\,} }

\title[Hessian Nilpotent Polynomials and Jacobian Conjecture]
{Hessian Nilpotent Polynomials and the Jacobian Conjecture}

    \author{Wenhua Zhao}      
    \date{}
    \begin{document}

\begin{abstract}

Let $z=(z_1, \cdots, z_n)$ and 
$\Delta=\sum_{i=1}^n \fr {\p^2}{\p z^2_i}$
the Laplace operator.
The main goal of the paper is 
to show that the well-known Jacobian conjecture 
without any additional conditions is 
equivalent to the following what we call 
{\it vanishing conjecture}: for any 
homogeneous polynomial $P(z)$ of degree $d=4$, 
if $\Delta^m P^m(z)=0$ for all $m \geq 1$, then
$\Delta^m P^{m+1}(z)=0$ when $m>>0$, or equivalently,
$\Delta^m P^{m+1}(z)=0$ when 
$m> \fr 32 (3^{n-2}-1)$. It is also shown in this paper 
that the condition $\Delta^m P^m(z)=0$ ($m \geq 1$)
above is equivalent to the condition that $P(z)$ is 
Hessian nilpotent, 
i.e. the Hessian matrix 
$\Hes P(z)=(\fr {\p^2 P}{\p z_i\p z_j})$ 
is nilpotent.
The goal is achieved by using the recent 
breakthrough work of M. de Bondt, A. van den Essen
\cite{BE1} and various results obtained 
in this paper on Hessian nilpotent polynomials.
Some further results on Hessian nilpotent polynomials
and the vanishing conjecture above   
are also derived.

\end{abstract}

\keywords{Hessian nilpotent polynomials, Deformed inversion pairs, 
the Heat equation, harmonic polynomials, the Jacobian conjecture.}
   
\subjclass[2000]{33C55, 39B32, 14R15, 31B05.}

 \bibliographystyle{alpha}
    \maketitle

\tableofcontents

\renewcommand{\theequation}{\thesection.\arabic{equation}}
\renewcommand{\therema}{\thesection.\arabic{rema}}
\setcounter{equation}{0}
\setcounter{rema}{0}
\setcounter{section}{0}

\renewcommand{\theequation}{\thesection.\arabic{equation}}
\renewcommand{\therema}{\thesection.\arabic{rema}}
\setcounter{equation}{0}
\setcounter{rema}{0}

\section{\bf Introduction}\label{S1}

Let $z=(z_1, z_2, \cdots, z_n)$ and 
$F(z)=z-H(z)$  be a formal map from $\bC^n$ to 
$\bC^n$ with $o(H(z))\geq 2$
and $G(z)$ the formal inverse map of $F(z)$.
 The well-known Jacobian conjecture 
first proposed by Keller \cite{Ke} in 1939 claims that, 
{\it if $F(z)$ is a polynomial map 
with the Jacobian $j(F)(z)=1$, the inverse map $G(z)$ 
must also be a polynomial map}. Despite intense study 
from mathematicians in more than half a century, 
the conjecture is still wide open
even for the case $n=2$. 
In 1998, S. Smale \cite{S} 
included the Jacobian conjecture 
in his list of $18$ important mathematical problems 
for $21$st century.
For more history and known results on
the Jacobian conjecture, 
see \cite{BCW}, \cite{E} and references there.
Recently, M. de Bondt and A. van den Essen
\cite{BE1} (Also see G. Meng \cite{M})
have made a breakthrough 
on the Jacobian conjecture.
They reduced the Jacobian conjecture to 
polynomial maps $F(z)=z-H(z)$ with 
$H(z)=\nabla P(z)=(\fr{\p P}{\p z_1}, \fr{\p P}{\p z_2}, \cdots ,\fr{\p P}{\p z_n})$
for some polynomials 
$P(z)\in\bC[z]$. In this paper,
we will refer to this reduction as 
the {\it gradient reduction} and 
the condition $H(z)=\nabla P(z)$ for some 
$P(z)\in\bC[[z]]$ as 
the {\it gradient condition}.
Note that, by Poincar\'e lemma,  a formal map $F(z)=z-H(z)$ with
$o(H(z))\geq 2$ satisfies 
the gradient condition if and only if its Jacobian matrix $JF(z)$ 
is symmetric. Following the terminology in \cite{BE1}, 
we also call the formal maps 
satisfying the gradient condition 
{\it symmetric} formal maps. 

For further discussion, let us fix the
following notions.
A power series 
$P(z)\in \bC[[z]]$ 
is said to be HN ({\it Hessian nilpotent})
if its Hessian matrix
$\Hes P(z)=(\fr {\p^2 P}{\p z_i\p z_j})$ 
is nilpotent. Let $t$ be a formal 
parameter which commutes with $z$.
The {\it deformed inversion pair} $Q_t(z)$ of any $P(z)\in \bC[[z]]$ with 
$o(P(z))\geq 2$ is the unique power series $Q_t(z)\in \bC[[z, t]]$ 
with $o(Q_t(z))\geq 2$ 
such that the formal map
$G(z)=z+ t\nabla Q_t(z)$ 
is the inverse map of $F(z)=z-t\nabla P(z)$. 
Recently, G. Meng \cite{M} and 
D. Wright \cite{Wr1} have derived 
a tree expansion formula for 
the inverse map 
$G(z)$ of $F(z)=F_{t=1}(z)$. 
In \cite{Z2}, by studying 
the PDE satisfied by the deformed 
inversion pairs $Q_t(z)$, the author has derived 
a recurrent formula and a
binary rooted tree 
expansion formula for the deformed inversion pairs $Q_t(z)$. 
Furthermore, some close relationships 
among the deformed inversion pairs, the Legendre transform, 
the inviscid Burgers' equations 
and the Jacobian conjecture are 
also clarified in \cite{Z2}. 
For some other recent results on symmetric 
polynomial or formal maps, 
see \cite{BE1}--\cite{BE5}, 
\cite{EW},  \cite{M}, \cite{Wr1}, \cite{Wr2} and \cite{Z2}.

In this paper,  we will use some 
general results in \cite{Z2} to  
study HNS (Hessian nilpotent power series) $P(z)$ 
and their deformed inversion pairs $Q_t(z)$. 
Furthermore, by using the {\it gradient reduction}
in \cite{BE1} and various results derived 
in this paper on HN polynomials, we will show that 
the Jacobian conjecture is equivalent to 
what we call  {\it vanishing conjectures} 
of homogeneous HN polynomials. (See the discussion below.)
We first derive the PDE's satisfied by $Q_t(z)$, 
$\Delta^k Q_t^m$ $(k, m\geq 1)$ and 
$\exp(sQ_t(z))$ $(s\in\bC^\times)$, 
where $\Delta=\sum_{i=1}^n \fr {\p^2}{\p z^2_i}$
is the Laplace operator. In particular, 
we show in Theorem \ref{Heat} 
that $\exp(sQ_t(z))$ $(s\in\bC^\times)$
is the unique power series 
solution of the Cauchy problem of 
the Heat equation with the initial condition 
$\left. \exp(sQ_t(z)) \right |_{t=0}=\exp(sP(z))$.
We then derive a uniform formula (See Theorem \ref{T3.4}.) 
for the powers $Q_t^k(z)$ $(k\geq 1)$ of the deformed 
inversion pairs $Q_t(z)$ of HNS $P(z)$. We also
prove a general theorem, Theorem \ref{T4.1.1},  
on a relationship between 
 $\{ \text{Tr\,} \text{Hes}^m (P(z)) | m\geq 1 \}$
and $\{\Delta^m P^m(z) | m\geq 1\}$ for any power series $P(z)$. 
From this theorem, 
we show in Theorem \ref{Crit-1} that,  
for any formal power series $P(z)$, it is HN 
if and only if 
$\Delta^m P^m(z)=0$ for any $ m\geq 1$, 
or equivalently, 
$\Delta^m P^m(z)=0$ for any $1\leq  m\leq n$. 
Finally, we prove  
some identities, vanishing properties and 
isotropic properties of 
$\{\Delta^k P^m(z) | m, k\geq 0 \}$ 
for HNS or HNP's (Hessian nilpotent polynomials) $P(z)$.
Some close relationships of 
the deformed pairs $Q_t(z)$ of HNS or HNP's $P(z)$
with the Heat equation and the Jacobian conjecture 
are also clarified.
In particular, we show that
the Jacobian conjecture without any additional conditions 
is equivalent to the following
 {\it vanishing conjectures}:  
for any HNP $P(z)$ of degree $d=4$, 
$\Delta^m P^{m+1}(z)=0$ for $m>>0$, 
or more precisely, for all 
$m> \fr 32 (3^{n-2}-1)$.

One remark is that, due to the identity 
$\text{Tr\,} \Hes (P)=\Delta P$, 
any HNS $P(z)$ is automatically harmonic, 
i.e. $\Delta P(z)=0$. 
Note that harmonic polynomials 
(See \cite{ABR}, \cite{H} and \cite{T}.) 
are among the most 
classical objects in mathematics 
and have been very well studied. 
The classical study on harmonic polynomials
started from Legendre, Laplace, Jacobi 
in the late eighteen century. 
The modern generalizations 
of harmonic polynomials, namely,  spherical functions,
were first studied by Cartan and 
Weyl in the 1930's and later by Gelfand, Harish-Chandra, etc.  
It is quite surprising to see that, first,
 HNP's as a family of very special
harmonic polynomials are closely 
related with the notorious Jacobian conjecture. 
Secondly, it seems that HNP's have been overlooked and 
have not been studied until the
recent work of M. de Bondt, 
A. van den Essen \cite{BE1} and G. Meng \cite{M}. 
Besides the connections with 
the Jacobian conjecture discussed above, 
another interesting aspect of HNP's  
is their connection with the classical 
inviscid Burgers' equation in Diffusion theory 
and also the Heat equation.
Actually, the vanishing conjecture above is 
also equivalent to 
saying that the power series solutions of certain Cauchy problems 
of the inviscid Burgers' equation and  
the Heat equation must be polynomials and the exponentials 
of polynomials, respectively. 
(See discussion in Section $4$ in \cite{Z2} and 
Conjecture \ref{Conj-3.1} in this paper.) 
Considering the
connections of HNP's with the classical 
objects described above, we believe that 
HNP's deserve much more attentions from mathematicians. 

Considering the length of this paper, we give the following 
detailed arrangement description. 
In Section \ref{S2}, we first fix some notation 
and definitions that are needed 
throughout the rest of this paper. 
We then briefly recall certain results 
obtained in \cite{Z2} and prove some 
preliminary results 
including the PDE (See Corollary \ref{C2.2.5})
satisfied by $\Delta^k Q_t^m$ $(k, m\geq 1)$ for 
the deformed inversion pairs $Q_t(z)$ of 
any power series $P(z)$.
In Section \ref{S3}, for any HNS $P(z)$ and 
its deformed inversion pair $Q_t(z)$, 
we derive the PDE's satisfied by $Q_t(z)$ 
and $\exp(sQ_t(z))$ $(s\in\bC^\times)$, 
from which we derive with two different proofs  
a uniform formula Eq.\,(\ref{E-Q-k}) 
for $Q_t^k(z)$ $(k\geq 1)$.  
In Section \ref{S4}, 
we prove a general theorem, Theorem \ref{T4.1.1},  
on a relationship between 
$\{ \text{Tr\,} \text{Hes}^m (P(z)) | m\geq 1 \}$
and $\{\Delta^m P^m(z) | m\geq 1\}$ 
for the universal formal power series 
$P(z)$ with $o(P(z))\geq 2$. 
From this theorem, 
we deduce a criterion in Theorem \ref{Crit-1} 
for the Hessian nilpotency of a formal power series $P(z)$ 
in terms of certain vanishing properties of 
$\{\Delta^m P^m(z) | m\geq 1\}$. 
%
%
In Section \ref{S5.1}, by using 
a fundamental theorem of harmonic polynomials (See Theorem \ref{T6.1.2}), 
we derive a criterion in Proposition \ref{Crit-2} for Hessian nilpotency of 
homogeneous harmonic polynomials. 
In Section \ref{S5.2}, 
we give constructions 
for some HNP's and HNS.
In Section \ref{S6}, by using 
some of the main results in the 
previous sections, 
we prove more properties of 
HNS or HNP's $P(z)$. 
We prove in Proposition \ref{P5.1.1}
an identity 
and in Theorem \ref{T5.2.2}
an equivalence of certain vanishing properties 
of $\{\Delta^k P^m(z) | k, m\geq 1\}$ 
for HNS $P(z)\in \bC[[z]]$.
In Theorem \ref{T5.3.3}, we show some
isotropic properties for homogeneous HNP's.  
In Section \ref{S7}, we discuss some applications to 
the Jacobian conjecture. We formulate the 
{\it vanishing conjecture}, Conjecture \ref{VC}, 
for (not necessary homogeneous) HNP's and 
the {\it homogeneous vanishing conjecture}, Conjecture \ref{HVC}, 
for homogeneous HNP's. We show in Proposition \ref{P6.4} 
that both conjectures above are equivalent to 
the Jacobian conjecture.

Finally, some remarks on this paper are as follows. 
First, for convenience, we will fix $\bC$ as our base field.
But, all results, formulas as well as their proofs
(except the 1st proof of Theorem \ref{T3.4})
obtained in this paper hold or work equally well if one 
replace $\bC$ by any $\bQ$-algebra. 
Secondly, 
we will not restrict our study just on HNP's. 
Instead, we will formulate and prove results for 
HN formal power series whenever 
they hold in this general setting.
Thirdly, for any HNP's or locally convergent 
HNS $P(z)$, 
all formal power series 
involved in this paper
are locally convergent. 
This can be seen either from 
the fact that any local analytic 
map with non-zero Jacobian 
at the origin has a locally convergent inverse, or from 
the well-known Cauchy-Kowaleskaya theorem 
(See \cite{R}, for example.)

{\bf Acknowledgment}:
The author is very grateful to Professor Arno van den Essen 
who has carefully read through the first preprint of this paper 
and pointed out many misprints and several mistakes.
Great thanks also go to Professor David Wright 
for personal communications, especially for informing
the author some of his own recent results.
The author also would like to thank Professor Mohan Kumar 
for personal communications.

\renewcommand{\theequation}{\thesection.\arabic{equation}}
\renewcommand{\therema}{\thesection.\arabic{rema}}
\setcounter{equation}{0}
\setcounter{rema}{0}

\section{\bf Deformed Inversion Pairs of Formal Power Series}\label{S2}

In this section, 
we first fix some notation and definitions that are needed in this paper. 
We then briefly recall certain results 
obtained in \cite{Z2} and prove some preliminary results. 

\subsection{Notation and Conventions} 

Once and for all, we fix the following notation and  conventions.

\begin{enumerate}
\item We fix $n\geq 1$ and set 
$z=(z_1, z_2, \cdots, z_n)$. For any $\bQ$-algebra $k$,
we denote by $k[z]$ (resp. $k[[z]]$) the
polynomial algebra (resp. formal power series algebra) 
over $k$ in $z_i$ $(1\leq i\leq n)$.

\item 
For any $\bQ$-algebra $k$,
by a formal map $F(z)$ from $k^n$ to $k^n$, we simply mean 
$F(z)=(F_1(z), F_2(z), \cdots, F_n(z))$
with $F_i(z)\in k[[z]]$ $(1\leq i\leq n)$. 
We denote by $J(F)$ and $j(F)$ 
the  Jacobian matrix and the Jacobian of $F(z)$, 
respectively.

\item We denote by $\Delta$ the Laplace operator 
$\sum_{i=1}^n \frac {\p^2}{\p z_i^2}$. Note that, 
a polynomial or formal power
series $P(z)$ is said to be {\it harmonic} if $\Delta P=0$.

\item For any $k\geq 1$ and
$U(z)=(U_1(z), U_2(z), \cdots, U_k(z)) \in \bC[[z]]^{\times k}$, 
we set
\begin{align*}
o(U(z))=\min_{1\leq i\leq k} o(U_i(z))
\end{align*}
and, when $U(z)\in \bC[z]^{\times k}$,
\begin{align*}
\deg U(z)=\max_{1\leq i\leq k} \deg U_i(z).
\end{align*}
For any $U_t(z) \in \bC[t][[z]]^{\times k}
\text{ or } \bC[[z, t]]^{\times k}$ ($k\geq 1$) 
for some formal parameter $t$, the notation 
$o(U_t(z))$ and $\deg U_t(z)$ always
stand for the order and 
the degree of $U_t(z)$ with respect to $z$, respectively.

\item For any $P(z)\in \bC[[z]]$, we denote by $\nabla P(z)$ the gradient 
of $P(z)$, i.e. $\nabla P=(\frac {\p P}{\p z_1}, \frac {\p P}{\p z_2}, 
\cdots, \frac {\p P}{\p z_n})$. We denote by $\Hes (P)(z)$ 
the Hessian matrix of $P(z)$, 
i.e. $\Hes(P)(z)=(\frac {\p^2 P(z)}{\p z_i\p z_j})$. 

\item All $n$-vectors in this paper are supposed to be 
column vectors unless stated otherwise.
For any vector or matrix $U$, we denote by $U^t$ its transpose.
The standard $\bC$-bilinear form
of $n$-vectors is denoted by  $<\cdot, \cdot >$.
\end{enumerate}

\vskip5mm
The following lemma will be very useful in our later arguments.  

\begin{lemma}\label{L2.1.2}
For any $P(z)\in \bC[[z]]$ and $m\geq 1$, we have
\begin{align}\label{E2.1.2}
\Delta P^{m+1}(z)=(m+1) P^{m}\Delta P+ m(m+1) P^{m-1} <\nabla P, \nabla P>  
\end{align}
or 
\begin{align}\label{E2.1.3}
P^{m-1} <\nabla P, \nabla P>  =\frac 1{m(m+1)}\left (
\Delta P^{m+1}-(m+1) P^{m}\Delta P\right ).
\end{align}

Furthermore, when $P(z)$ is harmonic, we have
\BQ\label{E2.1.4}
P^{m-1}\Delta P^2=\frac 2{m(m+1)}\Delta P^{m+1}.
\EQ
\end{lemma}

 \pf Consider
\begin{align*}
\Delta P^{m+1} &=
\sum_{i=1}^n \frac {\p }{\p z_i} \frac {\p }{\p z_i} P^{m+1} \\
&=(m+1) \sum_{i=1}^n \frac {\p }{\p z_i} ( P^{m} \frac {\p P}{\p z_i} ) \\
&=
(m+1) \sum_{i=1}^n  P^{m} \frac {\p^2 P}{\p z_i^2}  +
m(m+1) \sum_{i=1}^n P^{m-1} \frac {\p P}{\p z_i}  \frac {\p P }{\p z_i}\\
&=
(m+1)  P^{m} \Delta P  +
m(m+1)<\nabla P, \nabla P>P^{m-1}.
\end{align*}
Hence,  we get Eq.\,(\ref{E2.1.2}) and (\ref{E2.1.3}). 

Now suppose that $P(z)$ is harmonic, i.e. $\Delta P=0$.
By Eq.\,(\ref{E2.1.2}) with $m=1$, we have
 $\Delta P^2=2<\nabla P, \nabla P>$ or 
$<\nabla P, \nabla P>=\fr 12 \Delta P^2$.
It is easy to see that, in this case,  
Eq.\,(\ref{E2.1.4}) follows directly from 
Eq.\,(\ref{E2.1.3}).
\epfv

\subsection{Deformed Inversion Pairs of Formal Power Series} 

For any $P(z)\in \bC[[z]]$ with $o(P(z))\geq 2$, we set 
$F(z)=z-\nabla P(z)$. 
It is shown in $\S 1.1$ in \cite{M} 
(Also see Lemma $3.1$ in \cite{Z2}.) 
that there is a unique 
$Q(z)\in \bC[[z]]$ with $o(Q(z))\geq 2$ such that the
formal inverse of $F(z)$ is given by $G(z)=z+\nabla Q(z)$. 
We call $Q(z)$ the {\it inversion pair} of $P(z)$. Furthermore, 
following the arguments in \cite{Z1} and \cite{Z2}, we also consider 
the deformation $F_t(z)=z-t\nabla P(z)$ of $F(z)$, where $t$ is 
a formal parameter which commutes with variables $z_i$ $(1\leq i\leq n)$. 
By Lemma $3.1$ in \cite{Z2}, 
we know that there exists a unique 
$Q_t(z)\in \bC[[z, t]]$ with $o(Q_t(z))\geq 2$
such that the formal inverse $G_t(z)$ of $F_t(z)$ 
is given by $G_t(z)=z+t\nabla Q_t(z)$. 
Note that, when $o(P(z))\geq 3$, we actually have
$Q_t(z)\in \bC[t][[z]]$ and $Q_{t=1}(z)=Q(z)$.
In general,  $tQ_t(z)$ is nothing but 
the inversion pair of $tP(z)$ over 
the $\bQ$-algebra $\bC[[t]]$.
Another way to look at the inversion pair 
is as follows. Set 
$U(z)=\frac 12 \sum_{i=1}^n z_i^2-P(z)$ and 
$V(z)=\frac 12 \sum_{i=1}^n z_i^2+Q(z)$, 
then $V(z)$ is exactly 
the Legendre transform (See \cite{Ar}, \cite{M} and \cite{Z2}.) 
of the formal power series $U(z)$.

\begin{defi}\label{DIP}
For any $P(z)\in \bC[[z]]$ with $o(P(z))\geq 2$, 
$Q_t(z)\in \bC[[z, t]]$ defined above is called
the {\it deformed inversion pair} of $P(z)$.
\end{defi}

Another important definition is the following.

\begin{defi}
For any $P(z) \in \bC[[z]]$, 
we say $P(z)$ is HN $(${\it Hessian nilpotent}$)$ if its Hessian  
matrix $\Hes (P)=\left (\frac {\p^2 P}{\p z_i\p z_j}\right )$ is nilpotent.
\end{defi}

\begin{rmk}\label{RK2.4}
Note that, $\text{Tr\,}\Hes (P)=\Delta P$ for any $P(z)\in \bC[[z]]$. 
Hence any HN formal power series is harmonic. But the converse is not true.
For some examples of HNP's and HNS,  see Subsection \ref{S5.2}. 
\end{rmk}  

Throughout the rest of this paper, 
for any formal power series 
$P(z)\in \bC[[z]]$, we will fix the notation $F(z)$, $F_t(z)$, 
$G(z)$, $G_t(z)$, $Q(z)$ and $Q_t(z)$ 
defined above unless stated otherwise. 
We will use the short words HNS and HNP 
for ``HN power series" and ``HN polynomial", respectively.
Furthermore, we also define a sequence of formal power 
series $\{Q_{[m]}(z) | m\geq 1\}$ by writing 
\begin{align}\label{def-Qm}
Q_t(z)=\sum_{m=1}^\infty   t^{m-1} Q_{[m]}(z).
\end{align}

\begin{lemma}\label{L2.2.3}
For any formal power series $P(z)$, we have 
\begin{enumerate}
\item[(a)]  
\begin{align}\label{E2.4}
(\Delta Q_t)(F_t) =\sum_{k=1}^\infty t^{k-1} \text{Tr\,} \text{Hes}^k(P).
\end{align}
\item[(b)]  $P(z)$ is HN if and only if $Q_t(z)$ is harmonic as a formal power series in $z$, 
and if and only if $Q_t(z)$ is HN as a formal power series in $z$.
\end{enumerate}
\end{lemma}
\pf $(a)$ Set $N_t(z)=\nabla Q_t(z)$. It is easy to check that
\begin{align}
JN_t(z)=\Hes Q_t(z),\\
\text{Tr\,} \Hes (Q_t) =\Delta Q_t.
\end{align}
Then Eq.\,(\ref{E2.4}) follows directly from the equations above
and Eq.\,$(2.4)$ in \cite{Z2}.

$(b)$ follows directly by applying  
Lemma $2.2$ in \cite{Z2} to the formal map $F(z)=z-H(z)$ 
with $H(z)=\nabla P(z)$.
\epfv

The following theorem which was first proved in the 
unpublished preprint \cite{Z1} and later \cite{Z2} 
(See Theorem $3.6$ and Proposition $3.7$ in \cite{Z2}.) 
will play a fundamental role in this paper.

\begin{theo}\label{T2.2.4} \cite{Z2} \quad  
For any $Q_t(z)\in \bC[[z, t]]$ with $o(Q_t(z))\geq 2$ and 
$P(z)\in \bC[[z]]$ with $o(P(z))\geq 2$, 
the following statements are equivalent.
\begin{enumerate}
\item $Q_t(z)$ is the deformed inversion pair of $P(z)$.
\item $Q_t(z)$ is the unique power series solution of 
the following Cauchy problem of PDE's.
\begin{align}\label{Cauchy-2}
\begin{cases}
&\frac {\p Q_t(z)}{\p t}=\frac 12 <\nabla Q_t, \nabla Q_t>,\\
&Q_{t=0}(z)=P(z).
\end{cases}
\end{align}
\end{enumerate}
Furthermore, we have the following recurrent formula. 
\begin{align}
Q_{[1]}(z)&=P(z),\\
Q_{[m]}(z)&=\frac 1{2(m-1)} 
\sum_{\substack {k, l\geq 1 \\ k+l=m}}
<\nabla Q_{[k]}(z), \nabla Q_{[l]}(z)>
\end{align}
for any $m\geq 2$.
\end{theo}

\begin{corol}\label{C2.2.5}
For any $k\geq 0$ and $m\geq 1$, we have
\begin{align}\label{E2.2.7}
\frac {\p}{\p t} \Delta^k Q_t^m(z)=
\frac 1{2(m+1)} \Delta^{k+1} Q_t^{m+1} -\frac 12 \Delta^k (Q_t^m\Delta Q_t). 
\end{align}
\end{corol}

\pf Since $\fr \p{\p t}$ and $\Delta^k$ $(k\geq 0)$ commute, 
by applying $\Delta^k$ $(k\geq 1)$ to Eq.\,(\ref{E2.2.7}) with $k=0$, 
we get Eq.\,(\ref{E2.2.7}) for any $k\geq 1$. Therefore we may assume $k=0$.

Consider
\begin{align*}
\frac {\p}{\p t}  Q_t^m(z)&=m Q_t^{m-1}\frac {\p Q_t}{\p t} \\
\intertext{Applying Eq.\,(\ref{Cauchy-2}) in Theorem \ref{T2.2.4}:}
&=\frac m2 Q_t^{m-1} <\nabla Q_t, \nabla Q_t>\\
\intertext{Applying Eq.\,(\ref{E2.1.3}) in Lemma \ref{L2.1.2} to $Q_t(z)$:}
&=\frac m2 \frac 1{m(m+1)} (\Delta Q_t^{m+1}-(m+1)Q_t^m\Delta Q_t)\\
&=\frac 1{2(m+1)}  \Delta Q_t^{m+1} -\fr 12 Q_t^{m}\Delta Q_t. 
\end{align*}
\epfv

\renewcommand{\theequation}{\thesection.\arabic{equation}}
\renewcommand{\therema}{\thesection.\arabic{rema}}
\setcounter{equation}{0}
\setcounter{rema}{0}

\section{\bf Deformed Inversion Pairs of HN Power Series}\label{S3}

In this section, we study deformed inversion pairs $Q_t(z)$ 
of HNS (Hessian Nilpotent Formal Power Series) $P(z)\in \bC[[z]]$.
We first derive the PDE's satisfied 
by $Q_t(z)$ and $\exp (sQ_t(z))$ $(s\in \bC^ \times)$.
We then discuss some relationships 
among deformed inversion pairs, the Heat equation and 
the Jacobian conjecture.  
Note that similar relationships among 
deformed inversion pairs of formal power series (not necessarily HN), 
the inviscid Burgers' equations and the Jacobian conjecture have been 
discussed in \cite{Z2}.
Finally, we derive with two different proofs
a uniform non-recurrent formula 
(See Eq.\,(\ref{E-Q-k}).) for $Q_t^k(z)$ $(k\geq 1)$.

\begin{theo}\label{T3.1}
For any $Q_t(z)\in \bC[[z, t]]$ with $o(Q_t(z))\geq 2$ and 
HNS $P(z)\in \bC[[z]]$ 
with $o(P(z))\geq 2$, the following statements are equivalent.
\begin{enumerate}
\item   
$Q_t(z)$ is the deformed inversion pair of $P(z)$.
\item   $Q_t(z)$ is the unique power series solution of 
the following Cauchy problem of PDE.
\begin{align}
\begin{cases}\label{Cauchy-3}
\frac {\p Q_t(z)}{\p t}&=\frac 14 \Delta Q_t^2, \\
Q_{t=0}(z)&=P(z).
\end{cases}
\end{align}
\end{enumerate}
Furthermore, we have the following recurrent formula. 
\begin{align}
Q_{[1]}(z)&=P(z),\label{EEEE3.2} \\
Q_{[m]}(z)&=\frac 1{4(m-1)} \Delta
\sum_{\substack {k, l\geq 1 \\ k+l=m}}
 Q_{[k]}(z) Q_{[l]}(z) \label{EEEE3.3}
\end{align}
for any $m\geq 2$.
\end{theo}

\pf First, by Lemma \ref{L2.2.3}, $(b)$,  we have $\Delta Q_t(z)=0$.  
Therefore, for any $m\geq 1$,  $Q_{[m]}(z)$ is harmonic. 
Secondly, for any harmonic formal power series
$U(z), $V(z)$ \in \bC[[z]]$, it is easy to check that we have 
\BQ
\Delta (UV)=2 <\nabla U, \nabla V>.
\EQ
By using the facts above, it is easy to see that 
the implication $(1) \Rightarrow (2)$ and also 
the recurrent formulas Eq.\,(\ref{EEEE3.2}), 
(\ref{EEEE3.3}) follow directly from 
Theorem \ref{T2.2.4}.

To see $(2)\Rightarrow (1)$, 
we denote by $\widetilde Q_t(z)$ 
the deformed inversion pair of $P(z)$. 
By the fact proved above, we know that 
$\widetilde Q_t(z)$ also satisfies 
Eq.\,(\ref{Cauchy-3}). 
Since the power series 
solution of the Cauchy problem 
Eq.\,(\ref{Cauchy-3}) is unique, 
which is given recursively by 
Eq.\,(\ref{EEEE3.2}) and 
(\ref{EEEE3.3}),  we have 
$\widetilde Q_t(z)=Q_t(z)$. Therefore
$(2)\Rightarrow (1)$ also holds.
\epfv

A relation of deformed inversion pairs $Q_t(z)$ of HNS 
$P(z)$ with the Heat equation is given by the following theorem.

\begin{theo}\label{Heat}
Let $P(z)\in \bC [[z]]$ be HN with $o(P(z)) \geq 2$
and $Q_t(z)$ its deformed inversion pair. 
For any non-zero $s\in \bC$, set
\BQn
U_{t, s}(z)=\exp(sQ_t(z))=\sum_{k=0}^\infty \frac {s^k Q_t^k(z)}{k!}.
\EQn
Then, $U_{t, s}(z)$ is the unique formal power 
series solution of the following 
Cauchy problem of the Heat equation.
\begin{align}\label{Cauchy-5}
\begin{cases}
\frac {\p U_{t,s}} {\p t}(z) &= \frac 1{2s} \Delta U_{t, s}(z),\\
U_{t=0, s}(z)&= \exp(s P(z)).
\end{cases}
\end{align}
\end{theo}

\pf The uniqueness can be proved by 
viewing $U_{t, s}(z)$ as a power series in $t$ 
with coefficients in $\bC[[z]]$ and showing that 
the coefficients of $t^k$ $(k\geq 1)$ are 
recurrently determined by 
the coefficient of $t^0$ 
which is $U_{t=0, s}(z)= \exp(s P(z))$.
We skip the details here. 
For a similar argument, see the proof of 
Proposition $2.5$ in \cite{Z2}. 
Note that, when $P(z)$ is locally convergent, 
$Q_t(z)$ and $\exp(sQ_t(z))$ $(s\in \bC^\times)$ 
are locally convergent.
Then the uniqueness in this case
also follows from the 
Cauchy-Kowaleskaya theorem (See \cite{R}) 
in PDE. 

Now we show that $U_{t, s}(z)$ satisfies Eq.\,(\ref{Cauchy-5}).
First note that, the initial condition in
Eq.\,(\ref{Cauchy-5})
follows immediately from the one in Eq.\,(\ref{Cauchy-2}).
Secondly,  by Lemma \ref{L2.2.3}, $(b)$, 
we have $\Delta Q_t(z)=0$.

Consider 
\begin{align}\label{E2.4.26}
\frac {\p U_{t, s} }{\p t}&= \frac {\p \exp(s Q_t)}{\p t}\\
&=  s \exp(s Q_t) \frac {\p Q_t}{\p t} \nno \\
\intertext{Applying  Eq.\,(\ref{Cauchy-2}):}
&=  \frac s2 U_{t, s}  <\nabla  Q_t,  \nabla  Q_t>. \nno
\end{align}
On the other hand, we have 
\begin{align}\label{E2.4.27}
\Delta U_{t, s}&=\sum_{i=1}^n \frac{\p}{\p z_i}
\frac{\p}{\p z_i} e^{sQ_t} \\
&=s \sum_{i=1}^n \frac{\p}{\p z_i} (\frac{\p Q_t}{\p z_i} e^{sQ_t}) \nno\\
\intertext{Using the fact that $\Delta Q_t=0$:}
&=s^2 \sum_{i=1}^n \frac{\p Q_t}{\p z_i} \frac{\p Q_t}{\p z_i} e^{sQ_t}\nno \\
&=s^2  U_{s, t} <\nabla  Q_t,  \nabla  Q_t>. \nno
\end{align} 

By combining Eq.\,(\ref{E2.4.26}) and 
(\ref{E2.4.27}), we see that
$U_{t, s}(z)$ does satisfy the PDE in Eq.\,(\ref{Cauchy-5}). 
\epfv

By combining the gradient reduction in \cite{BE1}, \cite{M} 
and the homogeneous reduction in \cite{BCW},  \cite{Y} 
on the Jacobian conjecture, 
we see that the Jacobian conjecture 
can be reduced to polynomial maps $F(z)=z-\nabla P(z)$ 
with $P(z)\in \bC[z]$ homogeneous of degree $d=4$. 
By Theorem \ref{T3.1} and Theorem \ref{Heat}, 
it is easy to see that the Jacobian conjecture 
is equivalent to the following conjecture.

\begin{conj}\label{Conj-3.1}
For any homogeneous HNP $P(z)$ of degree $d\geq 2$, 
the unique solutions of the Cauchy problems 
$Eq.\,(\ref{Cauchy-5})$ and  
$Eq.\,(\ref{Cauchy-2})$ must be a polynomial in $(z, t)$ and 
the exponential of a polynomial $(z, t)$, respectively.
\end{conj}

Since it has been proved by Wang \cite{Wa} that 
the Jacobian conjecture holds for 
polynomial maps $F(z)$ of $\deg F(z)\leq 2$, 
hence Conjecture \ref{Conj-3.1} is true for $d\leq 3$. 
For more discussion on relationships of HNP's 
and the Jacobian conjecture, see Section \ref{S7}.

Next we give two different proofs for the following 
uniform formula for the powers $Q_t^k(z)$ $(k\geq 1)$
of the deformed inversion pairs $Q_t(z)$ of HNS $P(z)$. 

\begin{theo}\label{T3.4}
Suppose $P(z)\in \bC[[z]]$ with $o(P(z))\geq 2$ is HN. Then, 
for any $k \geq 1$, we have
\begin{align}\label{E-Q-k}
Q_t^k(z)=k! \sum_{m=0}^\infty \frac {t^m} {2^{m}m!(m+k)!} \Delta^{m} P^{m+k}(z).
\end{align}
In particular, for any $m \geq 1$, 
\begin{align}\label{E-Qm}
Q_{[m]}(z)= \frac 1{2^{m-1}m!(m-1)!} \Delta^{m-1} P^m(z).
\end{align}
\end{theo}

\underline{\it First Proof}: 
First, note that Eq.\,(\ref{E-Qm}) follows directly
from 
Eq.\,(\ref{E-Q-k}) with $k=1$
and the definition Eq.\,(\ref{def-Qm}) of $Q_{[m]}(z)$ $(m\geq 1)$. 
To prove Eq.\,(\ref{E-Q-k}), we consider 
the formal power series
\BQ\label{s-t}
\exp(\frac t{2s} \Delta) \exp (sP)=\sum_{k=0}^\infty \frac {t^k}{(2s)^k k!} 
\Delta^k \exp(sP).
\EQ
It is easy to check that the series above is also a 
formal power series solution of the Cauchy problem Eq.\,(\ref{Cauchy-5}). 
Hence, by Theorem \ref{Heat} and 
the uniqueness of the power series solution of 
Eq.\,(\ref{Cauchy-5}), we have
\BQ \label{E2.4.29}
\exp(sQ_t)=\exp(\frac t{2s} \Delta) \exp (sP).
\EQ
By comparing the coefficients of $s^k$ ($k\geq 1$) 
of the both sides of the equation above, 
we get Eq.\,(\ref{E-Q-k}).
\epfv

The proof above for Theorem \ref{T3.4} is shorter  
but less intriguing than the second proof below, which 
begins with the following lemma.
 
\begin{lemma}\label{L3.4}
Let $P(z)\in \bC[[z]]$ with $o(P(z))\geq 2$ be HN and 
$Q_t(z)$ the deformed inversion pair of $P(z)$.
Then, for any $k, l\geq 1$,  we have
\begin{align}\label{E2.4.31}
\frac {\p^{l} Q_t^k}{\p t^l}=
\frac {\Delta^l Q_t^{k+l}}{2^l (k+1)(k+2)\cdots (k+l)}. 
\end{align}
\end{lemma}

\pf We fix $k\geq 1$ and use the mathematical induction 
on $l\geq 1$. 
First, by Lemma \ref{L2.2.3}, $(b)$, 
we have $\Delta Q_t(z)=0$.  Then  
Eq.\,(\ref{E2.4.31}) for $l=1$ 
follows directly from 
Eq.\,(\ref{E2.2.7}). 

Now we assume that Eq.\,(\ref{E2.4.31}) 
holds for $l=l_0\geq 1$ and
consider the case $l=l_0+1$.

\begin{align*}
\frac {\p^{l} Q_t^k}{\p t^l} 
&=\frac {\p} {\p t}\left (\frac {\p^{l_0} Q_t^k}{\p t^{l_0}} \right )\\
&=\frac 1{2^{l_0} (k+1)(k+2)\cdots (k+l_0)} 
\frac {\p} {\p t} \Delta^{l_0} Q_t^{k+l_0} \\
&=\frac {k+l_0} {2^{l_0} (k+1)(k+2)\cdots (k+l_0)} 
 \Delta^{l_0}(Q_t^{k+l_0-1} \frac {\p Q_t} {\p t}) \\
\intertext{Applying the PDE in Eq.\,(\ref{Cauchy-3}):}
&=\frac  1{2^{l_0} (k+1)(k+2)\cdots (k+l_0-1)}\frac 14 
 \Delta^{l_0}(Q_t^{k+l_0-1} \Delta  Q_t^2 )\\
\intertext{Applying Eq.\,(\ref{E2.1.4}) to $Q_t^{k+l_0-1} \Delta  Q_t^2$:}
&=\frac  1{2^{l_0+2} (k+1) \cdots (k+l_0-1)}
\frac 2{(k+l_0)(k+l_0+1)}
 \Delta^{l_0+1}Q_t^{k+l_0+1} \\
&=\frac  1{2^{l_0+1} (k+1)(k+2)\cdots (k+l_0+1)}
 \Delta^{l_0+1}Q_t^{k+l_0+1}\\
&=\frac  1{2^{l} (k+1)(k+2)\cdots (k+l)}
 \Delta^{l}Q_t^{k+l}. 
\end{align*}
Hence, Eq.\,(\ref{E2.4.31}) holds for $l=l_0+1$. 
\epfv



\vskip3mm

\underline {\it 2nd Proof of Theorem \ref{T3.4}:}
First, by the initial condition in Eq.\,(\ref{Cauchy-2}), we have 
\begin{align}
\left. \Delta^l Q^{k+l}_{t}(z)\right |_{t=0}=\Delta^l P^{k+l}(z)
\end{align} 
for any $k\geq 1$ and $l\geq 0$.

Secondly, by setting $t=0$ in Eq.\,(\ref{E2.4.31}) and applying the equation above, 
we see that the coefficient of $t^l$ of $Q_t^k(z)$ is equal to 
\BQn
\frac 1{l!} \frac {\Delta^l P^{k+l}}{2^l (k+1)(k+2)\cdots (k+l)}
=\frac {k!}{2^l l!  (k+l)!} \Delta^l P^{k+l}.
\EQn
Hence Eq.\,(\ref{E-Q-k}) holds. 
\epfv

By comparing the coefficients of $s^k$ ($k\leq 0$) 
of the both sides of Eq.\,(\ref{s-t}), we see that
$\Delta^k P^m=0$ for any $k\geq m$, which is equivalent to saying that
$\Delta^m P^m=0$ for any $m\geq 1$. Note that the later statement  
also follows from Eq.\,(\ref{E-Q-k}) with $k=1$ and the fact 
that $Q_t(z)$ is harmonic.

\begin{corol} \label{C3.3}
For any HNS $P(z)\in \bC[[z]]$ with $o(P(z))\geq 2$, 
we have $\Delta^m P^m=0$ for any $m\geq 1$.
\end{corol}

Later, we will show in Theorem \ref{Crit-1} 
that the converse of the corollary above is also true.

Note that, by setting $s=1$ in Eq.\,(\ref{E2.4.29}), we have
the following formula.
\BQ
\exp (Q_t)=\exp ({\frac t2 \Delta})\exp (P).
\EQ

Actually, a more delicate formula (See Eq.\,(\ref{Exp-Form}) below.) 
can be derived as follows.

Set 
\begin{align}
&\widetilde Q_t(z)=Q_t(z)-P(z),\\
&\widetilde U_t(z)=\exp (\widetilde Q_t(z)), \\
&\Lambda_P(z) =\sum_{i=1}^n \frac {\p P(z)}{\p z_i}\frac {\p}{\p z_i}.
\end{align}
 
\begin{lemma}\label{L3.6}
For any HNS $P(z)\in \bC[[z]]$ with $o(P(z))\geq 2$, 
let $\widetilde Q_t(z)$, $\widetilde U_t(z)$ and $\Lambda_P$
be as above. Then we have

$(a)$ $\widetilde Q_t(z)$ is the unique power series solution of 
the following Cauchy problem of PDE's.
\begin{align}
\begin{cases}
\frac{\p \widetilde Q_t(z)}{\p t} 
&= \frac 14 \Delta \widetilde Q_t^2(z)+\Lambda_P \widetilde Q_t(z)+\frac 14 \Delta P^2,\\   
\widetilde Q_{t=0}(z)&=0.
\end{cases}
\end{align}

$(b)$ $\widetilde U_t(z)$ is the unique power series 
solution of the following Cauchy problem of PDE's.
\begin{align}\label{E2.4.35}
\begin{cases}
\frac{\p \widetilde U_t(z)}{\p t} 
&= (\frac 12 \Delta +\Lambda_P +\frac 14 \Delta P^2) \, \, \widetilde U_t(z),\\   
\widetilde U_{t=0}(z)&=1.
\end{cases}
\end{align}
\end{lemma}

The proof of this lemma is straightforward and similar as 
the proof of Theorem \ref{Heat}, so we omit it here. 
From Eq.\,(\ref{E2.4.35}), it is also easy 
to derive the following formula.

\begin{propo}\label{P3.7}
For any HNS $P(z)\in \bC[[z]]$ with $o(P(z))\geq 2$, we have
\begin{align}\label{Exp-Form}
\exp(\widetilde Q_t)=\exp (t(\frac 12 \Delta +\Lambda_P +\frac 14 \Delta P^2)) \cdot 1.
\end{align}
\end{propo}
\pf We first set 
\begin{align}
V_t(z)=\exp (t(\frac 12 \Delta +\Lambda_P +\frac 14 \Delta P^2)) \cdot 1.
\end{align}
 
Note that $V_{t=0}(z)=1$. Now we consider 
\begin{align*}
\fr {\p V_t}{\p t}&=(\frac 12 \Delta +\Lambda_P +\frac 14 \Delta P^2)
\exp (t(\frac 12 \Delta +\Lambda_P +\frac 14 \Delta P^2)) \cdot 1\\
&=(\frac 12 \Delta +\Lambda_P +\frac 14 \Delta P^2)\,\, V_t(z).
\end{align*}
Therefore $V_t(z)$ satisfies the PDE in 
Eq.\,(\ref{E2.4.35}). On the other hand, 
by Lemma \ref{L3.6}, $(b)$ above, $\exp(\widetilde Q_t)$ 
also satisfies Eq.\,(\ref{E2.4.35}). 
Hence, Eq.\,(\ref{Exp-Form}) follows from the uniqueness of
power series solutions of 
the Cauchy problem Eq.\,(\ref{E2.4.35}).
\epfv

One interesting aspect of the formula above is as follows. 
It shows the differential operator $\Lambda_P$ 
and the operator of the multiplication 
by $\Delta P^2$ also play important roles 
for deformed inversion pairs $Q_t(z)$. 
For example, from Eq.\,(\ref{Exp-Form}), 
it is easy to see that we have the following
corollary.

\begin{corol}\label{C3.9}
For any $P(z)\in \bC[[z]]$ with $o(P(z))\geq 2$ 
such that $\Delta^2 P(z)=0$, we have
 $\widetilde Q_t(z)=0$, 
or in other words, $Q_t(z)=P(z)$.
\end{corol}

\renewcommand{\theequation}{\thesection.\arabic{equation}}
\renewcommand{\therema}{\thesection.\arabic{rema}}
\setcounter{equation}{0}
\setcounter{rema}{0}

\section{\bf A Criterion for Hessian nilpotency} \label{S4}

Let $a=\{ a_I |I\in \bN^n, \, |I|\geq 2 \}$ be 
a set of variables which commute 
with each other. 
Let $P(z)=\sum_{I\in \bN^n } a_I z^I$ be the 
universal formal power series in $z$ with
$o(P(z))\geq 2$. We will also view $P(z)$ 
as a formal power series in $z$ 
with coefficients in $\bC[a]$, i.e. $P(z)\in \bC[a][[z]]$. 

For any $m\geq 1$, we set 
\begin{align}
u_m(P)&=\text{Tr\,}  \text{Hes}^m(P), \label{Def-u} \\
v_m(P)&=\Delta^m P^m.\label{Def-v}
\end{align}

In this section, we 
prove a general theorem, Theorem \ref{T4.1.1},  
about a relation between $\{u_m (P)|m\geq 1\}$ 
and $\{v_m (P)|m\geq 1\}$. 
Consequently, we get  
a criterion for Hessian nilpotency of formal power series 
$P(z)\in \bC[[z]]$ 
in terms of certain vanishing properties 
of $\{\Delta^m P^m (z)|m\geq 1\}$. 



Let $P(z)\in \bC[a][[z]]$,  $\{u_m(P) | m\geq 1\}$ and 
$\{ v_m(P) | m \geq 1 \}$ be defined as above. 

For any $k\geq 1$, we define ${\mathcal U}_k(P)$ 
(resp. ${\mathcal V}_k(P)$) to
be the ideal in $\bC[a][[z]]$ generated by  
$\{u_m(P) |1\leq m\leq k\}$ 
(resp. $\{ v_m(P) | 1\leq m\leq k \}$)
and all their partial derivatives of any order. 
For convenience, we also set 
${\mathcal U}_0(P)={\mathcal V}_0(P)=0$. 

The first main result of this subsection is the following theorem.

\begin{theo}\label{T4.1.1}
For any $k\geq 1$,
${\mathcal U}_k(P)={\mathcal V}_k(P)$ as 
ideals in $\bC[a][[z]]$. 
\end{theo}

One immediate consequence of 
Theorem \ref{T4.1.1} is the following corollary.

For any $k\geq 1$, we define $\wtilde {\mathcal U }_k(P)$ 
(resp. $\wtilde {\mathcal V}_k(P)$) to
be the ideal in $\bC[a]$ generated by all coefficients 
of $\{u_m(P) |1\leq m\leq k\}$ 
(resp. $\{ v_m(P) | 1\leq m\leq k \}$).  

\begin{corol}\label{C4.1.1}
\footnote{
Professor David Wright 
\cite{Wr3} recently informed 
the author that he has obtained a different 
proof for this result.} 
For any $k\geq 1$,  $\wtilde {\mathcal U}_k(P)=\wtilde {\mathcal V}_k(P)$
as ideals in $\bC[a]$.
\end{corol}

From Theorem \ref{T4.1.1} or Corollary \ref{C4.1.1}, it is easy to see 
that we have
the following criteria for 
Hessian nilpotency of formal power series $P(z)\in \bC[[z]]$.

\begin{theo}\label{Crit-1}
For any $P(z)\in \bC[[z]]$ with $o(P(z))\geq 2$, 
the following statements are equivalent.
\begin{enumerate}

\item[(1)] $P(z)$ is HN.
\item[(2)] $\Delta^m P^m=0$ for any $m\geq 1$.
\item[(3)] $\Delta^m P^m=0$ for any $1\leq m\leq n$.
\end{enumerate}
\end{theo}

\pf  By Theorem \ref{T4.1.1} or Corollary \ref{C4.1.1}, 
we have, for any fixed $k\geq 1$, 
$\Delta^m P^m=0$  $(1\leq m \leq k)$ if and only if 
$\text{Tr}\, \text{Hes}^m(P)=0$ $(1\leq m \leq k)$. 
Then, the theorem follows directly 
from the following facts in linear algebra, 
namely, 
for any $n\times n$ matrix $A\in M_n(\bC)$, 
$A$ is nilpotent 
if and only if $\text{Tr\,} A^m=0$ $( m \geq 1)$, and 
if and only if $\text{Tr\,} A^m=0$ $(1\leq m \leq n)$.
\epfv

Besides the criteria in Theorem \ref{Crit-1}, 
we  believe the following one is also true.

\begin{conj}\label{Conj-4.3}
Let $P(z)\in \bC [[z]]$ with $o(P(z))\geq 2$. 
If $\Delta^m P^m(z)=0$ for $m>>0$, then $P(z)$ is HN.
\end{conj}

In the rest of this subsection, we fix the universal 
formal power series
$P(z)=\sum_{I\in \bN^n} a_I z^I$ with $o(P(z))\geq 2$ 
and give a proof for 
Theorem \ref{T4.1.1}. 
Note that all results proved 
in the previous sections also hold 
for formal power series 
over the $\bC$-algebra $\bC[a]$. 
In particular, they hold for our
universal formal power series 
$P(z)\in \bC[a][[z]]$. 

We begin with the following two lemmas.

\begin{lemma}
Let $Q_t(z)\in \bC[a][[z, t]]$ be 
the deformed inversion pair of $P(z)$. 
Then, there exists a sequence $\{ w_k(P)(z) \in {\mathcal U}_k(P) |k \geq 1 \}$ 
such that
\begin{align}
\Delta Q_t(z)=\sum_{k=1}^\infty w_k(P)(z)\, \, t^{k-1}, \label{New-4.4}
\end{align}
and, for any $k\geq 1$,
\begin{align}
w_k(P)(z) \equiv u_k(P)(z)  \mod ({\mathcal U}_{k-1}(P)). \label{E2}
\end{align}
\end{lemma}

\pf 
We set $N_t(z)=\nabla Q_t(z)$.
First,  by composing $G_t(z)=z+tN_t(z)$ 
from right to 
Eq.\,(\ref{E2.4}) in Lemma \ref{L2.2.3}, 
we have 
\begin{align}\label{E1}
\Delta Q_t(z)=\sum_{i=1}^\infty u_i(P)(z+tN_t(z)) \, \, t^{i-1}.
\end{align}

Now we write the Taylor expansion of
$u_i(P)(z+tN_t(z))$ $(i\geq 1)$ at $z$ as
\begin{align}\label{Taylor}
u_i(P)(z+tN_t(z))=u_i(P)(z)+
\sum_{k=1}^\infty \sum_{\substack{ {\bf s}\in \bN^n \\ |{\bf s}|=k }}
\fr 1{{\bf s}!}\fr{\p^{\bf s}u_i(P)}{\p z^{\bf s}}(z)N_t^{\bf s}(z) t^k. 
\end{align}

Next we want to write the RHS of Eq.\,(\ref{Taylor}) 
as a formal power series in $t$ with coefficients in $\bC[[a, z]]$.
This can be done by first doing so for 
$N_t^{\bf s}(z)$ $({\bf s}\in \bN^n)$ in Eq.\,(\ref{Taylor})
and then re-arranging properly all the terms involved.
Note that, $u_i(P)(z)$ or $\frac{\p^{\bf s}u_i(P)}{\p z^{\bf s}}(z)$ 
$({\bf s }\in \bN^n)$
do not depend on $t$ and are in the ideal ${\mathcal U}_i(P)$. 
Also note that all the terms or products in the sum of 
Eq.\,(\ref{Taylor})
except the first one $u_i(P)(z)$
have positive degree in $t$ due to the factors 
$t^k$ $(k\geq 1)$. 
By using the observations above and 
keeping track the degree in $t$, 
it is easy to see that 
$u_i(P)(z+tN_t(z))$ can be written as

\begin{align}\label{E5}
u_i(P)(z+tN_t(z))=u_i(P)(z)+\sum_{j=1}^\infty A_{ij}(z)\, t^j
\end{align}
for some $A_{ij}(z)\in \mathcal U_i(P)(z)$ $(i, j\geq 1)$.

Now, by combining Eq.\,(\ref{E1}) and (\ref{E5}), we have

\begin{align*}
\Delta Q_t(z)&=\sum_{i=1}^\infty \left (u_i(P)(z)+\sum_{j=1}^\infty A_{ij}(z) t^j\right ) t^{i-1}\\
&=\sum_{k=1}^\infty (u_k(P)+\sum_{\substack{i+j=k\\i, j\geq 1}} A_{ij}(z) ) \, t^{k-1}\\
&=\sum_{i=1}^\infty w_k(P)(z)\,\, t^{k-1},  
\end{align*}
where, for any $k\geq 1$,
\begin{align}\label{E4.8}
w_k(P)(z)=u_k(P)+\sum_{\substack{i+j=k\\i, j\geq 1}} A_{ij}(z).
\end{align}
Hence we get Eq.\,(\ref{New-4.4}). 
By the fact that $A_{i, j}(z)\in \mathcal U_i(P)(z)$ for any $i, j\geq 1$, 
we see that each $A_{i, j}(z)$ in Eq.\,(\ref{E4.8}) lies in  
$\mathcal U_{k-1}(P)(z)$ since $i\leq k-1$. 
Therefore Eq.\,(\ref{E2}) also holds. 
\epfv

\begin{lemma}\label{L4.3}
For any $m, k\geq 1$ and $1\leq l\leq k$, 
we have
\begin{align}\label{E4.1.2}
 \frac {\p^l}{\p t^l} \Delta^m Q_t^m(z)
 \equiv \frac {\Delta^{m+l} Q_t^{m+l}}{2^l(m+1)\cdots (m+l)}
 \mod  ( {\mathcal U}_k(P), t^{k-l+1}).
\end{align}
\end{lemma}
\pf 
We fix $k\geq 1$ and use the mathematical induction on $l$ 
to show Eq.\,(\ref{E4.1.2}) holds for any $m\geq 1$. 

By Eq.\,(\ref{E2.2.7}) with $k=m$ and Eq.\,(\ref{New-4.4}), 
it is easy to see that
Eq.\,(\ref{E4.1.2}) holds for any $m\geq 1$ when $l=1$.

Now we assume that Eq.\,(\ref{E4.1.2}) holds 
for any $1\leq l\leq k_0<k$ and 
consider the case $l=k_0+1$. 
By applying $\frac {\p}{\p t}$ 
to Eq.\,(\ref{E4.1.2}) with $l=k_0$, we have
\begin{align}\label{E4.1.3}
  \frac {\p^{k_0+1}\Delta^m Q_t^m} {\p t^{k_0+1}}   
 \equiv &
 \frac {\frac {\p}{\p t} \Delta^{m+k_0} Q_t^{m+k_0}}{2^{k_0}(m+1)\cdots (m+k_0)}
\mod({\mathcal U}_k(P), t^{k-k_0}). 
\end{align}
While, from Eq.\,(\ref{E4.1.2}) with $l=1$, we have 

\begin{align}
\frac {\p}{\p t} \Delta^{m+k_0} Q_t^{m+k_0}\equiv
 \frac {\Delta^{m+k_0+1} Q_t^{m+k_0+1}}{2(m+k_0+1)}  \mod ({\mathcal U}_k(P), t^{k}).\nno
\end{align}

Since $k-k_0\leq k$,  hence we also have 

\begin{align}\label{E4.1.5}
 \frac {\p}{\p t} \Delta^{m+k_0} Q_t^{m+k_0}  \equiv  
\frac{\Delta^{m+k_0+1} Q_t^{m+k_0+1}}{2(m+k_0+1)}
 \mod ({\mathcal U}_k(P), t^{k-k_0}). 
\end{align}

By combining Eq.\,(\ref{E4.1.3}) and (\ref{E4.1.5}), we have 

\begin{align*}
 \frac {\p^{k_0+1}}
{\p t^{k_0+1}} \Delta^m Q_t^m  & \equiv  
 \frac {\Delta^{m+k_0+1} Q_t^{m+k_0+1}}   
{2^{k_0+1}(m+1)(m+2)\cdots (m+k_0+1)} \\
  & \quad \quad \quad \quad \quad \quad 
\mod ({\mathcal U}_k(P), t^{k-k_0}=t^{k-(k_0+1)+1})
\end{align*}
which is Eq.\,(\ref{E4.1.2}) for $l=k_0+1$.
\epfv

Now we are ready to prove Theorem \ref{T4.1.1}.
\vskip3mm

\underline{\it Proof of Theorem \ref{T4.1.1}:}
We use the mathematical induction on $k\geq 1$. Since 
$u_1(P)=\Delta P=v_1(P)$, hence, the theorem is 
 true for $k=1$. 

Now, we assume that ${\mathcal U}_k(P)={\mathcal V}_k(P)$  
for some $k\geq 1$. 
By Eq.\,(\ref{E4.1.2}) with $m=1$ and $l=k$, 
we have
\begin{align*}
\frac {\p^{k}}{\p t^{k}} \Delta Q_t \equiv 
\fr {\Delta^{k+1} Q_t^{k+1}}{2^{k}(k+1)!}
 \, \mod ({\mathcal U}_k(P), t).
\end{align*}
In other words, we have 
\begin{align}\label{E4.7}
\left. \frac {\p^{k}}{\p t^{k}} \Delta Q_t \right |_{t=0} 
\equiv \left. \fr {\Delta^{k+1} Q_t^{k+1}}{2^{k}(k+1)!}
\right |_{t=0}  \, \mod ({\mathcal U}_k(P)).
\end{align}
On the other hand, by  
the initial condition in Eq.\,(\ref{Cauchy-2}), we have
\begin{align*}
\left. \fr {\Delta^{k+1} Q_t^{k+1}}{2^{k}(k+1)!}
\right |_{t=0} &=\frac{\Delta^{k+1}P^{k+1}}{2^k(k+1)!}=\fr {v_{k+1}(P)}{2^k(k+1)!}.
\end{align*}
By  Eq.\,(\ref{New-4.4}) and  (\ref{E2}), we have
\begin{align*}
\left. \frac {\p^{k}}{\p t^{k}} \Delta Q_t \right |_{t=0} 
&=k!\, w_{k+1}(P)\equiv k! u_{k+1}(P) \mod ({\mathcal U}_k(P)).
\end{align*}
Therefore, by Eq.\,(\ref{E4.7}) and the two equations above, we have
\begin{align*}
u_{k+1}(P)&\equiv \fr 1{2^{k}k!(k+1)!} v_{k+1}(P) \, 
\mod ({\mathcal U}_k(P)).
\end{align*}
Since ${\mathcal U}_k(P)={\mathcal V}_k(P)$,
hence we have ${\mathcal U}_{k+1}(P)={\mathcal V}_{k+1} (P)$.
\epfv

\renewcommand{\theequation}{\thesection.\arabic{equation}}
\renewcommand{\therema}{\thesection.\arabic{rema}}
\setcounter{equation}{0}
\setcounter{rema}{0}

\section{\bf Hessian Nilpotent Polynomials}\label{S5}

In this section, we first derive in Subsection \ref{S5.1} 
a criterion for 
Hessian nilpotency of homogeneous polynomials 
by using a fundamental theorem 
(See Theorem \ref{T6.1.2}) 
of harmonic polynomials.
We then give in Subsection \ref{S5.2} 
some examples of  HNP's 
(Hessian nilpotent polynomials) and HNS 
(Hessian nilpotent formal power series). 


\subsection {A Criterion for Hessian Nilpotency of 
Homogeneous Harmonic Polynomials} \label{S5.1}

For any $n\geq 1$, 
we let $X(\bC^n)$ or simply $X$ denote the affine variety 
defined by $\sum_{i=1}^n z_i^2=0$. 
For any $d\geq 0$, 
we denoted by $V_d(z)$ the vector space of homogeneous polynomials in $z$ 
of degree $d\geq 0$.
For any $\alpha\in \bC^n$, 
we denote by $h_\alpha(z)$ 
the linear function 
$<\alpha, z>$ of $\bC^n$.

The following identities are almost trivial but very useful 
for our later arguments. 
So we formulate them as a lemma without giving proofs.

\begin{lemma}\label{L6.1.1}
$(a)$ For any $\alpha \in \bC^n$ and $m\geq 1$, we have
\begin{align}\label{E6.1.1}
\Hes (h_\alpha^m) (z)=  m(m-1) h_\alpha^{m-2} (z)\, \alpha \cdot \alpha^t
\end{align}

$(b)$
For any $\alpha, \beta \in X(\bC^n)$ and $m, k \geq 1$, we have
\begin{align}\label{E6.1.2}
\Delta (h_\alpha^m (z)h_\beta^k (z) )=2mk<\alpha, \beta >h_\alpha^{m-1} (z)h_\beta^{k-1} (z)
\end{align}
\end{lemma}

By Remark \ref{RK2.4}, we know that any HNS $P(z)$
is automatically harmonic, i.e. $\Delta P(z)=0$. 
For harmonic polynomials, we have 
the following fundamental theorem.

\begin{theo}\label{T6.1.2}
For any homogeneous harmonic polynomial 
$P(z)$ of degree $d\geq 2$, we have
\begin{align}\label{E6.1} 
P(z)=\sum_{i=1}^k c_i h_{\alpha_i}^d (z)
\end{align}
for some $c_i\in \bC$ and $\alpha_i\in X(\bC^n)$ $(1\leq i\leq k)$.
\end{theo}

For the proof of this theorem, see, for example, \cite{H} and \cite{T}.

Note that, by replacing 
$\alpha_i$ by $c_i^{-\fr 1d}\alpha_i$ $(1\leq i\leq k)$ 
in Eq.\,(\ref{E6.1}), 
we see that any  homogeneous harmonic polynomial 
$P(z)$ of degree $d\geq 2$ can be written as

\begin{align}
P(z)=\sum_{i=1}^k  h_{\alpha_i}^d (z) \label{d-Form}
\end{align}
for some $\alpha_i\in X(\bC^n)$ $(1\leq i\leq k)$.

In the rest of this subsection,  
we fix a homogeneous harmonic polynomial 
$P(z)\in V_d(z)$ of degree $d\geq 2$ and assume that 
$P(z)$ is given by Eq.\,(\ref{d-Form}) 
for some $\alpha_i \in X(\bC^n)$ $(1\leq i\leq k)$. 
We also assume   
$\{h_\alpha^d(z) | \alpha\in \bC^n \}$
are linearly independent in $V_d(z)$.

We first define the following matrices 
associated with $P(z)$.
\begin{align}
 A_P&=(<\alpha_i, \alpha_j>)_{k\times k}, \\
 \Psi_P&=(<\alpha_i, \alpha_j > h_{\alpha_j}^{d-2}(z) )_{k\times k}.
\end{align}

The main result of this section 
is the following proposition.

\begin{propo}\label{Crit-2}
Let $P(z)\in V_d(z)$ be given by Eq.\,$(\ref{d-Form})$. 
Then, for any $m\geq 1$, 
we have 
\begin{align}\label{E6.1.7}
\text{Tr\,} \Hes^m (P)=(d(d-1))^m \text{Tr\,} \Psi_P^m.
\end{align}
In particular, 
$P(z)$ is HN if and only if 
the matrix $\Psi_P$ is nilpotent.
\end{propo}
\pf
First, by Eq (\ref{d-Form}) and (\ref{E6.1.1}), 
we can write 
$\Hes (P)$ explicitly as 
\BQ\label{E6.2.8}
\Hes (P)=d(d-1) \sum_{i=1}^k 
 h_{\alpha_i}^{d-2}(z)\,\, \alpha_i \cdot \alpha_i^t.
\EQ

For any $m\geq 1$, we set $c_m=(d(d-1))^m$. By Eq.\,(\ref{E6.2.8}), 
we have

\begin{align*}
 \Hes^m (P) & =c_m\sum_{i_1, i_2, \cdots, i_m=1}^k 
 h_{\alpha_{i_1}}^{d-2}(z)
\cdots 
h_{\alpha_{i_m}}^{d-2}(z)
(\alpha_{i_1}  \cdot \alpha_{i_1}^t)
\cdots 
(\alpha_{i_m}  \cdot \alpha_{i_m}^t) \nno \\
& =c_m\sum_{i_1, i_2, \cdots, i_m=1}^k 
 \alpha_{i_1}(\alpha_{i_1}^t\cdot \alpha_{i_2}) (\alpha_{i_2}^t \cdot \alpha_{i_3})
\cdots (\alpha_{i_{m-1}}^t\cdot \alpha_{i_m})\alpha_{i_m}^t \nno\\
& \quad \quad   \quad \quad \quad \quad  
\quad \quad \quad   \quad \quad \quad \quad 
\cdot h_{\alpha_{i_1}}^{d-2}(z)
\cdots 
h_{\alpha_{i_m}}^{d-2}(z) \nno\\
 & =c_m\sum_{i_1, i_2, \cdots, i_m=1}^k 
< \alpha_{i_1}, \alpha_{i_2}> <\alpha_{i_2}, \alpha_{i_3}>
\cdots <\alpha_{i_{m-1}}, \alpha_{i_m}> \nno\\
& \quad \quad   \quad \quad \quad \quad  
\quad \quad \quad   \quad \quad \quad \quad 
\cdot h_{\alpha_{i_1}}^{d-2}(z)
\cdots 
h_{\alpha_{i_m}}^{d-2}(z)\,\,
\alpha_{i_1}\cdot \alpha_{i_m}^t \nno\\
& =c_m\sum_{i_1,  i_m=1}^k 
\left (\Psi_P^{m-1}\right )_{i_1, i_m}
h_{\alpha_{i_1}}^{d-2}(z)\,\,
\alpha_{i_1}\cdot \alpha_{i_m}^t. \nno
\end{align*}

By taking the trace of the matrices above, we get

\begin{align*}
 \text{Tr\,} \Hes ^m (P) 
&=c_m\sum_{i_1,  i_m=1}^k 
\left (\Psi_P^{m-1}\right )_{i_1, i_m}
\text{Tr\,} (\alpha_{i_1}\cdot \alpha_{i_m}^t)\, h_{\alpha_{i_1}}^{d-2}(z) \\
&=c_m\sum_{i_1,  i_m=1}^k 
\left (\Psi_P^{m-1}\right )_{i_1, i_m}
<\alpha_{i_m}, \alpha_{i_1}>  h_{\alpha_{i_1}}^{d-2}(z) \\
&=c_m\sum_{i_1,  i_m=1}^k 
\left (\Psi_P^{m-1}\right )_{i_1, i_m}
(\Psi_P)_{i_m, i_1}\\
&=c_m\text{Tr\,} \Psi_P^m.
\end{align*}
Hence, we get Eq.\,(\ref{E6.1.7}).
\epfv

\begin{corol}
Let $P(z)\in V_d(z)$ be given by Eq.\,$(\ref{d-Form})$. 
Suppose that $P(z)$ is HN. Then the matrix $A_P$ must be singular.
\end{corol}
\pf By Proposition \ref{Crit-2}, we have $\Psi_P$ is nilpotent. Therefore, we have
\begin{align}
0=\det \Psi_P
=h^{d-2}_{\alpha_1}(z) h^{d-2}_{\alpha_2}(z) \cdots h^{d-2}_{\alpha_k}(z)\det A_P.\nno  
\end{align}
Hence, we have $\det A_P=0$.
\epfv

\begin{corol}\label{C6.1.4}
Let $P(z)$ be HN and given by 
Eq.\,$(\ref{d-Form})$. Then, for any $2\leq m\leq d$, we have
\begin{align}\label{E6.1.9}
\sum_{i, j=1}^k <\alpha_i, \alpha_j>^m h_{\alpha_i}^{d-m}(z) h_{\alpha_j}^{d-m}(z)=0.
\end{align}
In particular, we have
\begin{align}\label{P-alpha=0}
\sum_{i=1}^k P(\alpha_i)=0.
\end{align}
\end{corol}

\pf First, note that, 
Eq.\,(\ref{P-alpha=0}) follows directly from  Eq.\,(\ref{E6.1.9}) with $m=d$.

To prove Eq.\,(\ref{E6.1.9}), we first consider the case $m=2$, 
the LHS of Eq.\,(\ref{E6.1.9}) is just $\text{Tr\,} \Psi_P^2$ 
up to a non-zero constant. 
Hence, by Proposition \ref{Crit-2}, Eq.\,(\ref{E6.1.9}) holds in this case.

Now consider the case $m>2$. By Eq.\,(\ref{E6.1.1}),
 we have
\begin{align}\label{E6.1.8}
 \Delta^l (h_{\alpha_i}^{d-2}(z) h_{\alpha_j}^{d-2}(z)) 
= 2^l & (d-2)^2 \cdots (d-l-1)^2 \cdot \\ 
& \quad \quad \cdot <\alpha_i, \alpha_j>^l 
 h_{\alpha_i}^{d-2-l}(z)h_{\alpha_j}^{d-2-l}(z) \nno
\end{align}
for any $1\leq i, j\leq k$ and $1\leq l\leq d-2$.

Now by applying $\Delta^{m-2}$ to 
Eq.\,(\ref{E6.1.9}) for the case $m=2$ and applying Eq.\,(\ref{E6.1.8}), 
it is easy to see that Eq.\,(\ref{E6.1.9}) holds for any $3\leq m\leq d-2$.
\epfv

Actually, by Eq.\,(\ref{E6.1.2}), 
the LHS of Eq.\,(\ref{E6.1.9}) is also $\Delta^m P^2(z)$ up to 
a non-zero constant. 
Hence the corollary above also follows from Theorem \ref{Crit-1}, 
which implies that $\Delta^m P^2(z)=0$ for any $m\geq 2$. 

One remark is that, by applying similar arguments 
as above to the equations $\Delta^m P^m(z)=0$ 
and $\text{Tr\,} \Psi_P^m=0$ ($m\geq 1$), 
one can derive more explicit identities 
satisfied by certain powers of $h_{\alpha_i}(z)$ 
$(1\leq i\leq k)$. 
But, in order to keep this paper 
in certain size, we skip them here. 
More study on homogeneous HNP's 
will be given in \cite{Z3}.

\subsection{Some Examples of HNP's and HNS} \label{S5.2}

In this subsection, we give some examples of HNS and HNP's. 

First, let $\Xi=\{\beta_i |1\leq i\leq k\}$
be any non-empty subset
of $X(\bC^n)$ such that $<\beta_i, \beta_j>=0$ 
for any $1\leq i, j\leq k$. 
For any $d \geq 2$, we set

\begin{align}\label{E6.2.11}   
W_{[\Xi, \, d]}(z) =\sum_{i=1}^k h_{\beta_i}^d(z). 
\end{align}

For convenience, we also set $W_{[\Xi, \, d]}(z)=0$ for $\Xi=\emptyset$.

Now let
$\widetilde \Xi=(\Xi_1, \Xi_2, \cdots, \Xi_m, \cdots )$ 
be a sequence of finite subsets of $X(\bC^n)$ such that, 
for any $m_1,  m_2  \geq 1$ and any 
$\beta_i \in \Xi_{m_i}$ $(i=1, 2)$, 
we have $<\beta_1, \beta_2>=0$. 
 We set
\begin{align}\label{E6.2.13}
W_{\widetilde \Xi} (z)  = \sum_{m=1}^\infty 
W_{[\Xi_m, \, m+1]}(z)=\sum_{m=1}^\infty 
\sum_{\beta_{m_i} \in \Xi_m} h_{\beta_{m_i}}^{m+1}(z). 
\end{align}
 
A more general construction is as follows. 

Let $w=(w_1, w_2, \cdots, w_k)$ be a sequence of commutative variables
and $( \beta_1, \beta_2, \cdots, \beta_k)$ a sequence of elements
of $\bC^n$ with $<\beta_i, \beta_j>=0$ $(1\leq i, j\leq k)$. 
For any formal power series $g(w)\in \bC[[w]]$, 
we define $U_g(z)\in \bC[[z]]$ by 
\begin{align}\label{g(w)}
U_g(z) =g(h_{\beta_1}(z), h_{\beta_2}(z), \cdots, h_{\beta_k}(z) ). 
\end{align}

One special case of the construction above is as follows.
We introduce new commutative variables  
$u=(u_1, u_2, ... , u_n)$ and $v=(v_1, v_2, ... , v_n)$. 
For any $g(z)\in \bC[[z]]$, we set 
\begin{align} \label{pg(uv)}
P_g (u, v) =g(u_1+\sqrt{-1}v_1, u_2+\sqrt{-1}v_2, \cdots,  u_n+\sqrt{-1}v_n).
\end{align}

Note that, by setting $w_i=z_i$ and $\beta_i\in X(\bC^{2n})$ 
such that $h_{\alpha_i}(u, v)=u_i+\sqrt{-1}v_i$ $(1\leq i\leq n)$, 
we have $U_g(u, v)=P_g(u, v)$.

The following lemma is easy to check 
directly by using Lemma \ref{L6.1.1}.

\begin{lemma}\label{L6.2.6}    
For any $P(z)\in \bC[[z]]$ given by Eq.\,$(\ref{E6.2.13})$, 
Eq.\,$(\ref{g(w)})$, 
we have 
\begin{align}
\Delta P^m(z)=0.  \label{E6.2.15}
\end{align}
In particular, $P(z)$ is HN.
\end{lemma}

Note that, by choosing  $\widetilde \Xi$ in 
Eq.\,(\ref{E6.2.13}) and 
$g(w)\in \bC[[w]]$ in  Eq.\,(\ref{g(w)}) properly, 
we can construct many HNS, HNP's and homogeneous HNP's.
Unfortunately, all these HNS or HNP's $P(z)$
are of ``trivial type" in the sense 
that their deformed inversion pair 
$Q_t(z)=P(z)$. This can be easily seen from 
Eq.\,(\ref{E6.2.15}) and Corollary \ref{C3.9}.  
A family of non-trivial HNP's was 
given in \cite{BE1} 
which was constructed as follows. 
  
Let $H(z)=(H_1(z), H_2(z), \cdots, H_n(z))\in \bC[z]^{\times n}$. 
Let $u, v$ as defined before 
Eq.(\ref{pg(uv)}) and set
\begin{align*}
P_H(u, v)=\sum_{i=1}^n  v_i H_i (u_1+\sqrt{-1}v_1, u_2+
\sqrt{-1}v_2, \cdots, u_n+\sqrt{-1}v_n ).
\end{align*}

It was shown in Lemma $1.2$ in \cite{BE1} 
that $P_H(u, v)\in \bC[u, v]$ is 
HN if and only if $JH(z)$ is nilpotent. 

\renewcommand{\theequation}{\thesection.\arabic{equation}}
\renewcommand{\therema}{\thesection.\arabic{rema}}
\setcounter{equation}{0}
\setcounter{rema}{0}

\section{\bf More Properties of HN Polynomials}\label{S6}

In this section, we derive more properties of HNS (Hessian nilpotent 
formal power series) and HNP's (Hessian nilpotent polynomials).
We prove an identity in Proposition \ref{P5.1.1} and 
an equivalence of certain vanishing properties 
in Theorem \ref{T5.2.2} of $\{\Delta^k P^m(z) | k, m\geq 1\}$ 
for HNS $P(z)\in \bC[[z]]$. In Subsection \ref{S6.3}, 
we study an isotropic property of $\{\Delta^k P^m(z) | k, m\geq 1\}$ 
for homogeneous HNP's $P(z)$. 

\subsection{An Identity of HN Formal Power Series}\label{S6.1}

Let $P(z)\in \bC[[z]]$ be a HNS. 
For any $k\geq 0$ and $\alpha \geq 1$, we set
\begin{align}\label{E5.1}
u_{k, \alpha}(P)=\fr {\alpha !}{2^k k!(k+\alpha)!}\Delta^k P^{k+\alpha}.
\end{align} 

\begin{propo}\label{P5.1.1}
For any $\alpha, \beta\geq 1$ and $m \geq 0$, we have
\begin{align}\label{E5.2}
u_{m, \alpha+\beta}(P)=\sum_{\substack{k+l=m\\ k, l\geq 0}} 
u_{k, \alpha}(P)u_{l, \beta}(P).
\end{align}
More explicitly, we have
\begin{align}\label{E5.3}
\Delta^m P^{m+\alpha+\beta}=\binom {\alpha+\beta}\alpha^{-1}
\sum_{\substack{k+l=m\\ k, l\geq 0}} 
\binom mk \binom {m+\alpha+\beta}{k+\alpha}
(\Delta^k P^{k+\alpha})(\Delta^l P^{l+\beta}).
\end{align}
\end{propo}

\pf First, it is easy to see that Eq.\,(\ref{E5.3}) follows directly from
Eq.\,(\ref{E5.2}) and (\ref{E5.1}). So we only need prove Eq.\,(\ref{E5.2}).

By Eq.\,(\ref{E-Q-k}) and (\ref{E5.1}), we have 
\begin{align}\label{E5.4}
Q_t^{\gamma}(z)=
\sum_{m=0}^\infty u_{m, \gamma}(P)t^m
\end{align}
for any $\gamma\geq 1$. By comparing the coefficients of $t^m$ 
of both sides of the equation $Q_t^{\alpha+\beta}(z)=Q_t^\alpha(z)Q_t^\beta (z)$, 
we see that Eq.\,(\ref{E5.2}) holds. 
\epfv

\subsection{A Vanishing Property of HN Formal Power Series}\label{S6.2}

\begin{theo}\label{T5.2.2}
For any HNS $P(z)\in \bC[[z]]$, 
the following statements are equivalent.
\begin{enumerate}
\item For any $k\geq 1$, $\Delta^m P^{m+k}=0$ when $m>>0$.
\item There exists $k_0 \geq 1$, $\Delta^m P^{m+k_0}=0$ when $m>>0$.
\item $\Delta^m P^{m+1}=0$ when $m>>0$.
\end{enumerate}
\end{theo}

\pf $(1)\Rightarrow (2)$ is trivial. To show
 $(2)\Rightarrow (3)$, we assume that 
$\Delta^m P^{m+k_0}=0$ when $m>M_0$ for some $M_0\geq 1$. 
For any $m>M_0+k_0-1$, we have
\BQn
\Delta^m P^{m+1}=\Delta^{k_0-1}(\Delta^{m-k_0+1} P^{(m-k_0+1)+k_0})=0.
\EQn
Hence $(3)$ holds in this case.

Now we consider $(3)\Rightarrow (1)$. Since $P(z)$ is HN, 
Eq.\,(\ref{E-Q-k}) in Theorem \ref{T3.4} holds for any $k\geq 1$.
In particular,  $Q_t(z)$ is a polynomial in $t$ 
with coefficients in $\bC[[z]]$ by our assumption of 
$(3)$. Therefore, for any $k\geq 1$, 
 $Q^k_t(z)$ is also a polynomial in $t$ 
with coefficients in $\bC[[z]]$. 
By Eq.\,(\ref{E-Q-k}) again, we see that $(1)$ holds. 
\epfv

We believe that Theorem \ref{T6.2} is still true 
without the Hessian nilpotency condition. 
Actually, if Conjecture \ref{Conj-4.3} is true, 
it is certainly the case. 
More precisely, suppose that one of the statements, say $(3)$, 
of Theorem \ref{T6.2} holds for some $P(z)\in \bC[[z]]$. 
Then we have 
\BQn
\Delta^{m+1}P^{m+1}=\Delta (\Delta^m P^{m+1})=0
\EQn
when $m>>0$. If Conjecture \ref{Conj-4.3} is true, then $P(z)$ is HN.
Hence all other statements of Theorem \ref{T6.2} also hold.

Later we will show in Theorem \ref{T6.2} 
that the Jacobian conjecture is equivalent 
to saying that one of the statements 
in Theorem \ref{T5.2.2} holds 
for HNP's $P(z)$.

\subsection{Isotropic Properties of 
Homogeneous HN Polynomials}\label{S6.3}

For any $1\leq i\leq n$, we set $D_i=\fr {\p}{\p z_i}$ 
and $D=(D_1, D_2,\cdots, D_n)$.
We define a $\bC$-bilinear map 
$\{\cdot, \cdot\}: \bC[z]\times \bC[z]\to \bC[z]$ by setting
$\{f, g\}=f(D)g(z)$
for any $f, g\in \bC[z]$.
The $\bC$-bilinear map $\{\cdot, \cdot\}$ defined above 
is closely related with
the following commonly used 
Hermitian inner product of $\bC[z]$. 
See, for example,  \cite{ABR}, \cite{H} and \cite{KR}. 

\begin{align*}
(\cdot, \cdot): \bC[z]\times \bC[z] &\to \bC, \\
(f,\quad g)\quad  &\to (f(D)\bar g)(0),
\end{align*}
where,  $\bar g(z)=\sum_{I\in \bN^n} \bar a_{\bf s} z^{\bf s}$
if $g(z)=\sum_{I\in \bN^n} a_{\bf s} z^{\bf s}$.
In particular,  for any homogeneous polynomials 
$f, g\in \bC[z]$ of the same degree, 
we have $\{f, \bar g\}=(f, g)$. 

Actually, the Hermitian inner product $(\cdot, \cdot)$ 
plays an very important role in the study of 
classical harmonic polynomials (See  \cite{ABR} and \cite{H}.).
Due to the connection of 
$\{\cdot, \cdot\}$ with 
the Hermitian inner product 
$(\cdot, \cdot)$ described above,  
we refer to the properties of HNP's 
derived in this subsection  
as certain isotropic properties. 

The main result of this subsection is the following theorem.

\begin{theo}\label{T5.3.3}
Let $P(z)$ be a homogeneous HNP of degree $d\geq 3$ and 
${\mathcal I} (P)$ the ideal of $\bC[z]$
generated by $\sigma^2:=\sum_{i=1}^n z_i^2$ 
and $\fr {\p P}{\p z_i}$ $(1\leq i\leq n)$.
Then, for any $f(z)\in {\mathcal I} (P)$ and $m\geq 0$, we have 
\begin{align}
\{f, \Delta^mP^{m+1}\}=f(D)\Delta^mP^{m+1}=0.
\end{align}
\end{theo}

To prove this theorem, we first need the following two lemmas.

\begin{lemma}\label{L5.3.4}
For any homogeneous polynomial $f(z)$ of degree $k\geq 1$, we have 
\begin{align}\label{E5.3.7}
& \quad \sum_{i_1, i_2, \cdots, i_k=1}^n \fr {\p^k f(z)}{\p z_{i_1}\p z_{i_2}\cdots\p z_{i_k}} 
\fr {\p^k }{\p z_{i_1}\p z_{i_2}\cdots\p z_{i_k}} \\
& =\sum_{\substack{ {\bf s} \in \bN^n \\ |{\bf s}|=k}} 
\binom{k}{\bf s} \fr {\p^{k} f}{\p z^{\bf s}}
\fr {\p^{k} }{\p z^{\bf s}} \nno\\ 
&= k! f(D),\nno 
\end{align}
where $\binom {k}{\bf s}=\frac {k!}{s_1! s_2!\cdots s_n!}$ for any 
${\bf s}=(s_1, s_2, \cdots, s_n)\in \bN^n$ with $|{\bf s}|=k$.
\end{lemma}
\pf Since Eq.\,(\ref{E5.3.7}) is linear on $f(z)$, we may assume that $f(z)$ is a single monomial, say, 
$f(z)=z_1^{l_1}z_2^{l_2}\cdots z_n^{l_n}$ with $l_i\geq 0$ and $\sum_{i=1}^n l_i=k$. Now, we consider
\begin{align*}
& \sum_{i_1, i_2, \cdots, i_k=1}^n\left(  \fr {\p^k }{\p z_{i_1}\p z_{i_2}\cdots\p z_{i_k}} 
z_1^{l_1}z_2^{l_2}\cdots z_n^{l_n}\right )
\fr {\p^k }{\p z_{i_1}\p z_{i_2}\cdots\p z_{i_k}}\\
&=\left( \fr {k!}{l_1!l_2!\cdots l_n!} \fr {\p^k }{\p z_1^{l_1}\p z_2^{l_2}\cdots \p z_n^{l_n}} 
z_1^{l_1}z_2^{l_2}\cdots z_n^{l_n}\right ) \fr {\p^k }{\p z_1^{l_1}\p z_2^{l_2}\cdots \p z_n^{l_n}} \\
&=k!\left (\fr {\p }{\p z_1}\right )^{l_1}\left (\fr \p {\p z_2}\right )^{l_2}\cdots 
\left( \fr \p{\p z_n}\right )^{l_n}\\
&=k! f(D).
\end{align*}
\epfv

\begin{lemma}\label{L5.3.5.5}
For any $f(z), g(z)\in \bC[[z]]$ and $l\geq 1$, we have 
\begin{align}\label{E5.10.10}
\Delta^l (g f)&=\sum_{\substack{k_1+ k_2+ k_3=l\\
k_1, k_2, k_3 \geq 0}} 2^{k_2}\binom {l}{k_1, k_2, k_3}
\sum_{i_1, i_2, \cdots, i_{k_2}=1}^n \fr {\p^{k_2} \Delta^{k_1}g(z)}{\p z_{i_1}\cdots\p z_{i_{k_2}}} 
\fr {\p^{k_2} \Delta^{k_3} f(z)}{\p z_{i_1}\cdots\p z_{i_{k_2}}}\\
&=\sum_{\substack{k_1+ k_2+ k_3=l\\
k_1, k_2, k_3 \geq 0}} 2^{k_2}\binom {l}{k_1, k_2, k_3}
\sum_{\substack{ {\bf s} \in \bN^n \\ |{\bf s}|=k_2}} 
\binom {k_2}{\bf s} \fr {\p^{k_2} \Delta^{k_1} g}{\p z^{\bf s}}
\fr {\p^{k_2}\Delta^{k_3} f }{\p z^{\bf s}}.\nno
\end{align}
where $\binom {l}{k_1, k_2, k_3}=\frac {l!}{k_1! k_2! k_3!}$ for any $k_1, k_2, k_3\geq 0$. 
\end{lemma}

\pf We use the mathematical induction on $l\geq 1$. 
When $l=1$, by the Leibniz's rule, it is easy to check that 
\begin{align} \label{E5.10.10-1}
\Delta (gf)(z)=(\Delta g(z))f(z)+2\sum_{i=1}^n \fr{\p g(z)}{\p z_i} \fr{\p f(z)}{\p z_i}
+g(z)\Delta f(z)
\end{align}
which is exactly 
Eq.\,(\ref{E5.10.10}) with $l=1$.

Now we assume Eq.\,(\ref{E5.10.10}) holds
for $l=l_0\geq 0$. 
By using Eq.\,(\ref{E5.10.10-1}) and Eq.\,(\ref{E5.10.10}) with $l=l_0$, 
we have
\begin{align*}
& \quad \Delta^{l_0+1} (g f) = \Delta (\Delta^{l_0} (g f)) \\
& =  \sum_{\substack{k_1+ k_2+ k_3=l\\
k_1, k_2, k_3 \geq 0}} 2^{k_2}\binom {l}{k_1, k_2, k_3}
\sum_{i_1, i_2, \cdots, i_{k_2}=1}^n 
\Delta \left (
\fr {\p^{k_2} \Delta^{k_1}g(z)}{\p z_{i_1}\cdots\p z_{i_{k_2}}} 
\fr {\p^{k_2} \Delta^{k_3} f(z)}{\p z_{i_1}\cdots\p z_{i_{k_2}}} \right ) \nno \\
& =  \sum_{\substack{k_1+ k_2+ k_3=l\\
k_1, k_2, k_3 \geq 0}} 2^{k_2}\binom {l}{k_1, k_2, k_3}
\sum_{i_1, i_2, \cdots, i_{k_2}=1}^n 
\fr {\p^{k_2} \Delta^{k_1+1}g(z)}{\p z_{i_1}\cdots\p z_{i_{k_2}}} 
\fr {\p^{k_2} \Delta^{k_3} f(z)}{\p z_{i_1}\cdots\p z_{i_{k_2}}}  \nno \\
&  +  \sum_{\substack{k_1+ k_2+ k_3=l\\
k_1, k_2, k_3 \geq 0}} 2^{k_2+1}\binom {l}{k_1, k_2, k_3} 
\sum_{i_1, i_2, \cdots, i_{k_2+1}=1}^n 
\fr {\p^{k_2+1} \Delta^{k_1}g(z)}{\p z_{i_1}\cdots\p z_{i_{k_2+1}}} 
\fr {\p^{k_2+1} \Delta^{k_3} f(z)}{\p z_{i_1}\cdots\p z_{i_{k_2+1}}} \nno \\
&  +   \sum_{\substack{k_1+ k_2+ k_3=l\\
k_1, k_2, k_3 \geq 0}} 2^{k_2}\binom {l}{k_1, k_2, k_3}
\sum_{i_1, i_2, \cdots, i_{k_2}=1}^n 
\fr {\p^{k_2} \Delta^{k_1}g(z)}{\p z_{i_1}\cdots\p z_{i_{k_2}}} 
\fr {\p^{k_2} \Delta^{k_3+1} f(z)}{\p z_{i_1}\cdots\p z_{i_{k_2}}} \nno \\
\intertext{By shifting indices and using the convention that
$\binom {l}{k_1, k_2, k_3}=0$ if $k_1$, $k_2$ or $k_3<0$ :}
& =  \sum_{\substack{k_1+ k_2+ k_3=l+1\\
k_1, k_2, k_3 \geq 0}} 2^{k_2} \left ( \binom {l}{k_1-1, k_2, k_3}+
\binom {l}{k_1, k_2-1, k_3} \right. \nno \\
&  \quad \quad\quad\quad \left.
+\binom {l}{k_1, k_2, k_3-1} \right ) 
 \sum_{i_1, i_2, \cdots, i_{k_2}=1}^n 
\fr {\p^{k_2} \Delta^{k_1}g(z)}{\p z_{i_1}\cdots\p z_{i_{k_2}}} 
\fr {\p^{k_2} \Delta^{k_3} f(z)}{\p z_{i_1}\cdots\p z_{i_{k_2}}}.\nno
\end{align*}

Hence, we only need show that 
\begin{align*}
  \binom {l+1}{k_1, k_2, k_3}= \binom {l}{k_1-1, k_2, k_3}+
\binom {l}{k_1, k_2-1, k_3} 
+\binom {l}{k_1, k_2, k_3-1}
\end{align*}
for any $k_1, k_2, k_3\geq 0$. But this identity follows directly from the binomal 
expression of $(x+y+z)^k$ for $k=l, l+1$ and the identity 
\begin{align*}
(x+y+z)^{l+1}=(x+y+z)(x+y+z)^l.
\end{align*}
\epfv

One immediate consequence of Eq.\,$(\ref{E5.10.10})$ is the following corollary.
\begin{corol}
For $g(z), f(z)\in \bC[[z]]$ with $g(z)$ harmonic, 
we have
\begin{align} \label{E5.10}
\Delta^l (g f)=\sum_{k=0}^l 2^k \binom lk
\sum_{\substack{ {\bf s} \in \bN^n \\ |{\bf s}|=k}} \binom{k}{\bf s}
\fr {\p^{k} g}{\p z^{\bf s}}
\fr {\p^{k} \Delta^{l-k} f}{\p z^{\bf s}}. 
\end{align}
\end{corol}

Now we are ready to prove Theorem \ref{T5.3.3}.
\vskip3mm

\underline{\it Proof of Theorem \ref{T5.3.3}:}
Since $P(z)$ is HN, by Theorem \ref{Crit-1}, 
we know that $\Delta^m P^m=0$ for any $m\geq 1$. 
Therefore, we have, 
\begin{align*}
\sigma^2 (D) \Delta^m P^{m+1}
=\Delta^{m+1} P^{m+1}=0.
\end{align*}
Hence it will be enough 
to prove the theorem  
for $f(z)=\fr {\p P}{\p z_i}(z)$ 
 for any fixed $1\leq i\leq n$.

Since $d\geq 3$, we have $m+d-1 \geq  m+2$ and  
\begin{align}
\Delta^{m+d-1} P^{m+2}=0.
\end{align}
By applying $\fr \p{\p z_i}$ to the equation above,  
we get 
\begin{align}\label{EE5.3.1}
\Delta^{m+d-1} \fr {\p P}{\p z_i}P^{m+1} =\Delta^{m+d-1} ( f P^{m+1})=0.
\end{align}

But, on the other hand, by Eq.\,(\ref{E5.10}), we have
\begin{align} \label{EEE6.12}
& \quad \Delta^{m+d-1} \left ( f P^{m+1}\right )\\
& =
\sum_{k=0}^{m+d-1} 2^k \binom {m+d-1}k
\sum_{\substack{ {\bf s} \in \bN^n \\ |{\bf s}|=k}} \binom k{\bf s} 
\fr {\p^{k} f}{\p z^{\bf s}}
\fr {\p^{k} \Delta^{m+d-1-k} P^{m+1} }{\p z^{\bf s}}. \nno
\end{align} 
Furthermore, we also have the following equations.  
\begin{align*}
& \fr {\p^{\bf s} f}{\p z^{\bf s}}=0 \quad \quad 
\text{  for any ${\bf s} \in \bN^n$ with $|{\bf s}|= k > d-1$, } \\ 
& \Delta^a P^{m+1}=0 \quad  \quad \text{ for any $a\geq m+1$.}
\end{align*} 
Hence all the terms in the RHS of Eq.\,(\ref{EEE6.12}) 
except the one with $k=d-1$ are zero. Therefore,
by Eq.\,(\ref{E5.3.7}) in Lemma \ref{L5.3.4},  
we have
\begin{align*}
\Delta^{m+d-1} (f P^{m+1})
& =2^{d-1}\sum_{\substack{ {\bf s} \in \bN^n \\ |{\bf s}|=d-1}} 
\binom{d-1}{\bf s} \fr {\p^{d-1} f}{\p z^{\bf s}}
\fr {\p^{d-1} \Delta^{m} P^{m+1}}{\p z^{\bf s}} \\ 
& =2^{d-1} (d-1)! f(D) \Delta^{m} P^{m+1}.
\end{align*}
Hence, by Eq.\,(\ref{EE5.3.1}) and the equation above, 
we see that the theorem holds for 
$f(z)=\fr {\p P}{\p z_i}(z)$ $(1\leq i\leq n)$.
\epfv

\begin{corol}\label{C5.3.6}
Let $P(z)$ be a homogeneous HNP of degree $d\geq 3$ and 
$Q_t(z)$ its deformed inversion pair. 
Then, for any $k, l\geq 0$ with $k>l$, we have
\begin{align}\label{E5.3.12}
(\Delta^l P^k)(D)Q_t(z)=0. 
\end{align}
More precisely, we have
\begin{align}\label{E5.3.13}
(\Delta^l P^k)(D)  (\Delta^m P^{m+1}(z)) =0
\end{align} 
for any $k, l\geq 0$ with $k >l$.

In particular, we have
\begin{align}\label{E5.3.14}
P(D) Q_t(z)=0.
\end{align}
\end{corol}

\pf First note that Eq.\,(\ref{E5.3.14}) follows from 
Eq.\,(\ref{E5.3.12}) by setting $l=0$ and $k=1$. 
By Eq.\,(\ref{E-Q-k}), we see that Eq.\,(\ref{E5.3.12}) and
Eq.\,(\ref{E5.3.13}) are equivalent to each other. 
Hence, it is enough to show Eq.\,(\ref{E5.3.13}).
Furthermore, by Theorem \ref{T5.3.3}, 
it will be enough to show that $\Delta^l P^k(z)$ 
for any $k, l\geq 0$ with $k>l$
lies in the ideal $\tilde {\mathcal I}(P)$ 
generated by $\fr {\p P}{\p z_i}(z)$ $(1\leq i\leq n)$. 

By Euler's lemma, we have 
$P(z)=\fr 1d \sum_{i=1}^n z_i \fr{\p P}{\p z_i}(z)$. 
Hence, for any $k\geq 1$, 
$P^k(z)\in \tilde {\mathcal I}(P)$ and 
Eq.\,(\ref{E5.3.13}) holds when $l=0$.
Now we consider the case $l>0$. Note that
$\Delta^l$ is a sum 
of the differential operators of 
the form 
$\fr{\p^2}{\p z^2_{i_1}}\fr{\p^2}{\p z^2_{i_2}}\cdots \fr{\p^2}{\p z^2_{i_l}}$
with $1\leq i_1, i_2,\cdots, i_l\leq n$.
When we distribute the $2l$ derivations 
$\fr{\p}{\p z_{i_j}}$ $(1\leq j\leq l)$
of the differential operator above to $k$ 
copies $P(z)$ of $P^k(z)$, there is always at least 
one copy $P(z)$ of $P^k(z)$ receives 
none or one derivation. Otherwise, 
we would have $2k \leq 2l$, 
which contradicts to our condition $k>l$. 
Since we have already shown $P(z)\in \tilde {\mathcal I}(P)$
above, hence we have $\Delta^l P^k(z)\in \tilde {\mathcal I}(P)$ 
for any  $k, l\geq 0$ with $k>l>0$.
\epfv 
 
Theorem \ref{T5.3.3} and Corollary \ref{C5.3.6} do not hold 
for homogeneous HNP's $P(z)$ of degree $d=2$. But, by similar arguments 
as the proof of Theorem \ref{T5.3.3} starting from 
$\Delta^{m+2} P^{m+2}=\Delta^{m+2} (P\cdot P^{m+1})=0$
instead of Eq.\,(\ref{EE5.3.1}), one can show the following proposition.

\begin{propo}
Let $P(z)$ be a homogeneous HNP of degree $d=2$ and 
${\mathcal J} (P)$ the ideal of $\bC[z]$
generated by $P(z)$ and $\sigma^2=\sum_{i=1}^n z_i^2$.
Then, for any $f(z)\in {\mathcal J} (P)$ and $m\geq 0$, we have 
\begin{align}
f(D)\Delta^mP^{m+1}=0.
\end{align}
In particular, we have 
\begin{align}
P(D) Q_t(z)=0.
\end{align}
\end{propo}

\renewcommand{\theequation}{\thesection.\arabic{equation}}
\renewcommand{\therema}{\thesection.\arabic{rema}}
\setcounter{equation}{0}
\setcounter{rema}{0}

\section{\bf  Vanishing Conjectures of HN Polynomials}\label{S7}

In this section, we propose some conjectures on 
certain vanishing properties of polynomials 
$\{\Delta^k P^m(z) | m, k\geq 1 \text{ and } m>k \}$
for HNP's (Hessian Nilpotent polynomials) $P(z)$. 
We also show that these so-called {\it vanishing conjectures }
are equivalent to the well-known Jacobian conjecture.
 
\begin{conj} \label{VC} {\bf (Vanishing Conjecture )} \newline
For any HN $(\text{not necessarily homogeneous})$ 
polynomial $P(z)$ of degree $d\geq 2$, 
its deformed inversion pair $Q_t(z)$ is a polynomial in both $t$ and $z$. 
More precisely, $\Delta^k P^{k+1}=0$ when $k>>0$.
\end{conj}

\begin{theo}\label{T6.2}
The following statements are equivalent.
\begin{enumerate}
\item The vanishing conjecture for homogeneous HNP of degree $d=4$.
\item The vanishing conjecture for homogeneous HNP of degree $d\geq 2$.
\item The vanishing conjecture.
\item The Jacobian conjecture.
\end{enumerate}
\end{theo}
\pf First, it is easy to see that $(2)\Rightarrow (1)$,  
$(3)\Rightarrow (1)$ and $(3)\Rightarrow (2)$
are trivial. 
By the gradient reduction in \cite{BE1}
and the homogeneous reduction in \cite{BCW},  \cite{Y}
on the Jacobian conjecture, we know that the Jacobian conjecture will be true 
if it is true for polynomial maps $F(z)=z-\nabla P(z)$ with $P(z)$ being 
homogeneous HNP of degree $d=4$. Therefore, 
$(4)\Leftrightarrow (1)$ 
follows directly from Eq.\,(\ref{E-Q-k}) in Theorem \ref{T3.4}. 
Hence we only need show $(4)\Rightarrow (3)$.

Now we assume the Jacobian conjecture and let $P(z)$ be 
a HNP of degree $d\geq 2$. 
Let $F_t(z)=z-t\nabla P(z)$ and 
$G_t(z)=z+t\nabla Q_t(z)$ as before.
Consider the formal map 
$U(z, t)=(F_t(z), t)$ from $\bC^{n+1} \to \bC^{n+1}$.
It is easy to check that 
the Jacobian of the
map $U(z, t)$ with respect to $(z, t)$ 
is also identically equal to $1$
and the formal inverse $V(z, t)$ 
is given by $V(z, t)=(G_t(z), t)$. Since we have assumed 
the Jacobian conjecture, $V(z, t)$ must 
be a polynomial in $(z, t)$. Hence so is $G_t(z)$. 
By Eq.\,(\ref{E-Q-k}) again,   we see that 
$\Delta^m P^{m+1}(z)$ must vanish when $m>>0$. 
\epfv

Next, by using the upper bound 
given in \cite{BCW}, Corollary $1.4$, 
for the degrees of inverse maps of 
polynomial automorphisms of $\bC^n$, 
we show that Conjecture \ref{VC} for  
homogeneous HNP's can actually be 
reformulated more precisely as follows.

\begin{conj}\label{HVC}{\bf (Homogeneous Vanishing Conjecture)}\newline
For any homogeneous HNP $P(z)$ of degree $d\geq 2$, we have
\begin{enumerate}
\item \label{HVC1} $\Delta^m P^{m+1}=0$ for any $m> \alpha_{[n, d]}:=\fr 1{d-2}((d-1)^{n-1}-(d-1))$.
\item \label{HVC2} For any $k\geq 1$, $\Delta^m P^{m+k}=0$ for any $m> k\alpha_{[n, d]}$.
\end{enumerate}
\end{conj}

\begin{propo}\label{P6.4} 
$(a)$ For any homogeneous HNP $P(z)$ of 
degree $d\geq 2$, the statements $(1)$ and $(2)$ in Conjecture \ref{HVC} are equivalent.

$(b)$ Conjecture \ref{HVC} for $d\geq 2$ and Conjecture \ref{HVC} for $d=4$ 
are both equivalent to the Jacobian conjecture.
\end{propo}
\pf $(a)$ First, $(2)\Rightarrow (1)$ is trivial. Now we show that $(1)\Rightarrow (2)$. 
By Eq.\,(\ref{E-Q-k}) with $k=1$, we see that the degree $\deg_t Q_t(z)$ of $Q_t(z)$ with
respect to $t$ is less or equal to 
$\alpha_{n, d}$, i.e. $\deg_t Q_t(z)\leq \alpha_{[n, d]}$. Therefore, for any $k\geq 1$, 
we have,  $\deg_t Q^k_t(z)\leq k\alpha_{[n, d]}$. By Eq.\,(\ref{E-Q-k}) again, 
we see that $(2)$ holds for any $k\geq 1$.

$(b)$ By theorem \ref{T6.2}, it is easy to see that Conjecture \ref{HVC} for $d\geq 2$
or Conjecture \ref{HVC} with $d=4$ implies the Jacobian conjecture. 
Therefore it will be enough 
to show that the Jacobian conjecture implies Conjecture \ref{HVC}. 
By Corollary $1.4$ in \cite{BCW}, we know that, for any polynomial automorphism $F(z)$
of $\bC^n$ with $\deg F(z)=d-1$, $\deg G(z) \leq (d-1)^{n-1}$. 
 By applying this result to the polynomial automorphism 
$F(z)=z-\nabla P(z)$ and its
inverse $G(z)=z-\nabla Q(z)$,  
we get $\deg Q(z)\leq (d-1)^{n-1}+1$.  
Furthermore, 
by Eq.\,(\ref{E-Q-k}) with $k=1$ and the fact 
$\deg \Delta^m P^{m+1}(z)=(m+1)d-2m$  $(m\geq 1)$, we get 
$\Delta^m P^{m+1}(z)=0$ if
\begin{align*}
(m+1)d-2m >(d-1)^{n-1}+1.
\end{align*}
By separating $m$ from the inequality above, we get
$\Delta^m P^{m+1}(z)=0$ whenever $m>\alpha_{[n, d]}$.
\epfv

Finally, by translating certain known results 
on the Jacobian conjecture, 
we know that the vanishing conjectures, 
Conjecture \ref{VC} or \ref{HVC}, are true 
for the following cases.
{\it
\begin{enumerate}
\item[$\bullet$]  S. Wang \cite{Wa} proved that the Jacobian
conjecture holds for any polynomial map $F(z)$ of $\deg F(z)\leq 2$. 
Hence Conjecture \ref{VC} and \ref{HVC} hold for any HNP $P(z)$ of 
degree $d\leq 3$.

\item[$\bullet$] For any symmetric polynomial map 
$F(z)=z-H(z)$ with $o(H(z))\geq 2$ and $JH(z)$ nilpotent, 
A. van den Essen and S. Washburn \cite{EW} showed that 
the Jacobian conjecture holds when $n\leq 4$ and 
$H(z)$ is homogeneous. Later, M. de Bondt and 
A. van den Essen \cite{BE2}-\cite{BE4} further 
proved that the Jacobian conjecture holds 
either $n\leq 4$ without $H(z)$ 
being homogeneous, 
or $n=5$ with $H(z)$ being homogeneous. 
$($For an exposition discussion on these results, 
see \cite{BE5}.$)$. From the results above, we see that, 
Conjecture \ref{VC} has an affirmative answer
when $n\leq 4$, and Conjecture \ref{HVC} 
is true when $n\leq 5$. 

\item[$\bullet$] By Theorem $4.1$ in \cite{EW} 
and similar arguments there,  it is easy to show that 
the only HNP's $($not necessarily homogeneous$)$ 
$P(z)\in {\mathbb R}[z]$ 
$(o(P(z))\geq 2)$ with real coefficients
are $P(z)=0$. Hence, Conjecture \ref{VC} and \ref{HVC} 
hold trivially in this case.

\item[$\bullet$] Recently, D. Wright \cite{Wr1} showed that 
the Jacobian conjecture holds for any symmetric 
polynomial map $F(z)=z-H(z)$ with $H(z)$ homogeneous 
and $JH^3(z)=0$. Hence Conjecture \ref{HVC} holds
for any homogeneous HNP $P(z)$ with $\Hes ^3(P(z))=0$.
\end{enumerate}
}

Furthermore, by Corollary \ref{C3.9}, 
we also know that Conjecture \ref{VC} 
holds for any HNP $P(z)$ such that 
$\Delta P^2(z)=0$. In particular, it is true 
for HNP's constructed by Eq.(\ref{E6.2.13}) and Eq.(\ref{g(w)}).


{\small \sc Department of Mathematics, Illinois State University,
Normal, IL 61790-4520.}

{\em E-mail}: wzhao@ilstu.edu.

\end{document}